# ON A CLASS OF OPTIMAL STOPPING PROBLEMS FOR DIFFUSIONS WITH DISCONTINUOUS COEFFICIENTS

By Ludger Rüschendorf and Mikhail A. Urusov[1]

*University of Freiburg and Berlin University of Technology*

In this paper, we introduce a modification of the free boundary problem related to optimal stopping problems for diffusion processes. This modification allows the application of this PDE method in cases where the usual regularity assumptions on the coefficients and on the gain function are not satisfied. We apply this method to the optimal stopping of integral functionals with exponential discount of the form $E_x \int_0^\tau e^{-\lambda s} f(X_s)\,ds$, $\lambda \geq 0$ for one-dimensional diffusions $X$. We prove a general verification theorem which justifies the modified version of the free boundary problem. In the case of no drift and discount, the free boundary problem allows to give a complete and explicit discussion of the stopping problem.

**1. Introduction.** This paper is concerned with a class of optimal stopping problems for integral functionals with exponential discount $E_x \int_0^\tau e^{-\lambda s} f(X_s)\,ds$ for a one-dimensional diffusion process $(X_s)$. An example is the problem of stopping a Brownian motion as close as possible to its maximum as it can be reduced to this kind of stopping problem [see Graversen, Peskir and Shiryaev (2000) or Peskir and Shiryaev (2006), Chapter 8, Section 3.1]. The literature on optimal stopping problems for diffusion processes is very rich. An effective method for solving problems of this type is to develop a connection with some related free boundary problems (of Stefan type). There are two types of results on this connection. One type of result is based on solving the related free boundary problem with the smooth fit condition (in certain cases, the smooth fit condition is replaced by the continuous fit condition, also, additional conditions may be necessary). This allows to find explicit solutions of the initial optimal stopping problem in certain cases. Many examples of this approach (with explicit solutions) are presented in the recent comprehensive book of Peskir and Shiryaev (2006).

[1]This paper was written while the second author held a research fellowship of the Alexander von Humboldt Foundation.

Further, for some classes of optimal stopping problems for regular diffusion processes with smooth coefficients, existence and regularity results for the corresponding free boundary problems have been established under different kinds of smoothness and boundedness conditions on the coefficients of the diffusion and on the gain function.

On the other hand, starting with the work of Bensoussan and Lions (1973), general existence and regularity results for solutions of optimal stopping problems in terms of variational inequalities have been established. These are formulated with differential operators in the weak sense and allow weaker assumptions on the regularity of the coefficients and on the gain function. We refer to the book of Friedman (1976) for strong results in this direction [see also Nagai (1978) and Zabczyk (1984)]. Results of this type have led, in particular, to the development of some effective numerical solution methods [see, e.g., Glowinski, Lions and Trémolières (1976) and Zhang (1994)]. The method of variational inequalities is, however, typically more difficult to use in concrete examples, where one wants to find explicit solutions, compared to the formulation in terms of the free boundary PDE's.

In our paper, we discuss in detail the optimal stopping of integral functionals $E_x \int_0^\tau e^{-\lambda s} f(X_s)\, ds$ in the case of not necessarily continuous coefficients of the diffusion $X$ and for an interesting class of not necessarily continuous (cumulative) gain functions $f$. It is therefore obvious that the classical formulation of the free boundary problem is not applicable to these stopping problems. In Section 2, a suitably generalized formulation of the free boundary problem is given and a verification theorem, together with uniqueness results for this free boundary problem and for the optimal stopping time (Theorem 2.1), are proved. An important point in establishing this verification theorem is to establish variational inequalities for the solutions of our generalized free boundary problem (Lemma 2.6).

We would like to mention some related papers of Salminen (1985), Beibel and Lerche (2000), Dayanik and Karatzas (2003) and Dayanik (2003), where problems of maximizing $E_x[e^{-A_\tau} g(X_\tau) I(\tau < \infty)]$ over all stopping times $\tau$ are studied ($A$ is a continuous additive functional of $X$). These authors use different approaches, obtain some general results and explicitly treat several examples. Their approaches are also applicable to diffusions with discontinuous coefficients. Neither of these approaches is based on the free boundary method (the one we use here). Another difference with our paper is that we consider optimal stopping of *integral* functionals.

After finishing this paper, we became aware of the work Lamberton and Zervos (2006), where some interesting results about the value function of the problem of maximizing $E_x[e^{-A_\tau} g(X_\tau) I(\tau < \infty)]$ are proved.

Neither the function $g$ nor the coefficients of the diffusion $X$ are supposed to be continuous. Lamberton and Zervos (2006) prove that the value function $V$ in the problem they consider is the difference of two convex functions and satisfies a certain variational inequality. The results in our paper go in the opposite direction. Theorem 2.1 states that the solution of a certain (modified) free boundary problem is the value function $V$ in the problem we consider. One of the conditions in this free boundary problem is that $V$ should be differentiable and $V'$ absolutely continuous. The weaker condition mentioned above, that $V$ is the difference of two convex functions, is not sufficient for Theorem 2.1. The free boundary problem loses the uniqueness property and there appear solutions of the free boundary problem that have nothing to do with the stopping problem we consider [see Remark (i) after Theorem 2.1]. The reason is that the integral functionals $E_x \int_0^\tau e^{-\lambda s} f(X_s)\, ds$ are "more regular" than the functionals $E_x[e^{-A_\tau} g(X_\tau) I(\tau < \infty)]$. Hence, one may expect that value functions for integral functionals should also be "more regular."

In Section 3, we study in complete detail the case of diffusions without drift and with zero discount. Necessary and sufficient conditions for the existence of solutions of the free boundary problem are established (Theorem 3.1). Also, in the case that the free boundary problem has no solutions, the optimal stopping problem is dealt with (Theorems 3.2 and 3.9). We discuss finiteness of the value function, obtain explicit formulas for the value function and the optimal stopping time and determine approximately optimal stopping time sequences in the case where there is no optimal stopping time. Finally, the Appendix contains several technical lemmas which are used in the proofs and can be helpful in studying related questions.

## 2. Stopping problem, free boundary problem and verification theorem.

2.1. *Setting of the problem.* Let $X = (X_t)_{t\in[0,\infty)}$ be a continuous stochastic process with values in the extended real line $\mathbb{R} \cup \{\Delta\}$ and explosion time $\zeta$, that is, the following two properties hold:

(i) $X$ is $\mathbb{R}$-valued and continuous on $[0, \zeta)$;

(ii) if $\zeta < \infty$, then $X \equiv \Delta$ on $[\zeta, \infty)$ and either $\lim_{t\uparrow\zeta} X_t = \infty$ or $\lim_{t\uparrow\zeta} X_t = -\infty$.

If $\zeta = \infty$ a.s., then $X$ is called *nonexplosive.*

We consider stopping problems for diffusions $X$ defined by

$$dX_t = b(X_t)\, dt + \sigma(X_t)\, dB_t, \tag{2.1}$$

where $B$ is a Brownian motion and $b$, $\sigma$ are Borel functions $\mathbb{R} \to \mathbb{R}$. In the sequel, we assume that the coefficients $b$ and $\sigma$ satisfy the Engelbert–Schmidt

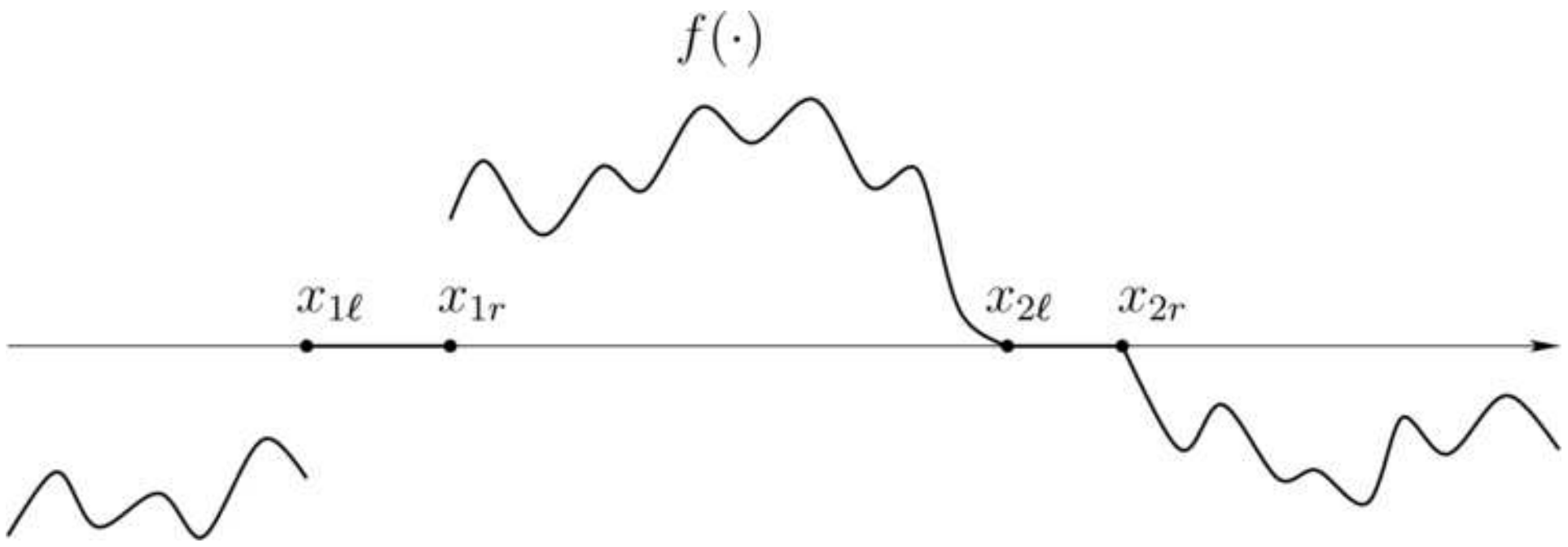

Fig. 1.

condition

$$\sigma(x) \neq 0 \ \ \forall x \in \mathbb{R} \qquad \frac{1}{\sigma^2} \in L^1_{\text{loc}}(\mathbb{R}), \qquad \frac{b}{\sigma^2} \in L^1_{\text{loc}}(\mathbb{R}), \tag{2.2}$$

where $L^1_{\text{loc}}(\mathbb{R})$ denotes the class of locally integrable functions on $\mathbb{R}$. Under this condition, the SDE (2.1) has a unique (in law, possibly explosive) weak solution for any starting point $X_0 = x$ [see Engelbert and Schmidt (1985, 1991) or Karatzas and Shreve (1991), Chapter 5, Theorem 5.15]. Let us note that condition (2.2) is weak enough. For example, if $b$ is locally bounded and $\sigma$ is locally bounded away from zero, then (2.2) holds.

Let $f:\mathbb{R} \to \mathbb{R}$ be a Borel function such that there exist points $x_{1\ell} \le x_{1r} < x_{2\ell} \le x_{2r}$ in $\mathbb{R}$ such that $f > 0$ on $(x_{1r}, x_{2\ell})$, $f = 0$ on $[x_{1\ell}, x_{1r}] \cup [x_{2\ell}, x_{2r}]$ and $f < 0$ on $(-\infty, x_{1\ell}) \cup (x_{2r}, \infty)$ (see Figure 1). Throughout the following, we assume that the function $f$ satisfies the condition

$$\frac{f}{\sigma^2} \in L^1_{\text{loc}}(\mathbb{R}), \qquad \frac{1}{f} \in B_{\text{loc}}(x) \qquad \forall x \in \mathbb{R} \setminus ([x_{1\ell}, x_{1r}] \cup [x_{2\ell}, x_{2r}]), \tag{2.3}$$

where $B_{\text{loc}}(x)$ denotes the class of functions $\mathbb{R} \to \mathbb{R}$ locally bounded at $x$ [$g \in B_{\text{loc}}(x)$ if it is bounded in a sufficiently small neighborhood of $x$]. Finally, let us note that $f/\sigma^2 \in L^1_{\text{loc}}(\mathbb{R})$ holds for all locally bounded functions $f$, due to (2.2).

In the paper, we use the following convention. For any function $g:\mathbb{R} \to \mathbb{R}$, we define $g(\Delta) = 0$.

Let $X$ be a (possibly explosive) solution of (2.1) on some probability space $(\Omega, \mathcal{F}, P_x)$, where $P_x(X_0 = x) = 1$, $x \in \mathbb{R}$. By $(\mathcal{F}^X_t)$, we denote the filtration generated by $X$ satisfying the usual conditions. In this paper, we consider the class of optimal stopping problems defined for functions $f$ as specified above by

$$V^*(x) = \sup_{\tau \in \mathfrak{M}} E_x \int_0^\tau e^{-\lambda s} f(X_s)\, ds. \tag{2.4}$$

Here, $\lambda \ge 0$ and $\mathfrak{M}$ is the class of $(\mathcal{F}_t^X)$-stopping times $\tau$ such that

$$E_x \int_0^\tau e^{-\lambda s} f^+(X_s)\,ds < \infty \quad \text{or} \quad E_x \int_0^\tau e^{-\lambda s} f^-(X_s)\,ds < \infty. \tag{2.5}$$

Let us remark that it is enough to consider only stopping times $\tau \le \zeta$ because $f(\Delta) = 0$. Setting (2.4) is a well-motivated class of stopping problems arising, for example, in connection with various versions of options of Asian type. The gain function $f$ positively rewards the case where the process stays in a favorable domain $[x_{1\ell}, x_{2r}]$ and puts negative weight in the case where this domain is left. From this formulation, one may expect that two-sided stopping times play an essential role.

For real numbers $\alpha < \beta$, we denote by $T_{\alpha,\beta}$ the stopping time

$$T_{\alpha,\beta} = \inf\{t \in [0,\infty) : X_t \le \alpha \text{ or } X_t \ge \beta\} \tag{2.6}$$

(as usual, $\inf \varnothing = \infty$). It is important that, under our assumptions, $T_{\alpha,\beta} \in \mathfrak{M}$ (see Lemma A.3 in the Appendix).

Let $\mu_L$ denote the Lebesgue measure on $(\mathbb{R}, \mathcal{B}(\mathbb{R}))$. Using the occupation times formula [see Revuz and Yor (1999), Chapter VI, Corollary (1.6)], one can verify that

$$\mu_L(\{t \in [0,\infty) : X_t \in A\}) = 0, \qquad P_x\text{-a.s.}$$

for sets $A$ of Lebesgue measure 0. Hence, problem (2.4) remains unchanged if we change the function $f$ on a set of $\mu_L$-measure 0. In particular, the cases, where some of the values $f(x_{1\ell})$, $f(x_{1r})$, $f(x_{2\ell})$ and $f(x_{2r})$ are nonzero, reduce to the situation under consideration.

The interesting point in the formulation of assumptions (2.2) and (2.3) is that they allow discontinuities in $b$, $\sigma$ and $f$, which is of interest in various applications. For example, in modeling stock prices, it is reasonable to consider volatilities $\sigma$ that jump to a higher range of values when the price reaches a certain threshold.

REMARKS. (i) All results of Section 2 also remain valid for $J \cup \{\Delta\}$-valued diffusions $X$, where $J = (\ell, r)$, $-\infty \le \ell < r \le \infty$, the functions $b$, $\sigma$ and $f : J \to \mathbb{R}$ satisfy conditions similar to (2.2) and (2.3) (one should replace $\mathbb{R}$ with $J$) and $X$ explodes when it tends either to $\ell$ or to $r$ at a finite time.

(ii) One main reason why we consider functionals of the form shown in Figure 1 is that one encounters concrete stopping problems of this type in the literature. For example, Graversen, Peskir and Shiryaev (2000) reduce the problem of stopping a Brownian motion $B = (B_t)_{t\in[0,1]}$ as close as possible to its maximum

$$\inf_\tau E\Big(\max_{t\in[0,1]} B_t - B_\tau\Big)^2$$

to the problem of the form (2.4) with an Ornstein–Uhlenbeck process $X$, $dX_t = X_t\,dt + \sqrt{2}\,dB_t$, $\lambda = 2$ and $f(x) = 3 - 4\Phi(|x|)$, where $\Phi$ is the distribution function of the standard Gaussian random variable [see also Peskir and Shiryaev (2006), Chapter VIII, Section 30.1]. Note that $f$ has the form shown in Figure 1.

As another example, we consider the following stopping problem of Karatzas and Ocone (2002), which they study in order to solve a stochastic control problem [see also Dayanik and Karatzas (2003), Section 6.9]:

$$\inf_{\tau} E_x\left[e^{-\lambda\tau}\delta X_\tau^2 + \int_0^\tau e^{-\lambda s} X_s^2\,ds\right], \qquad x \in J := (0,\infty), \tag{2.7}$$

where $\lambda \ge 0$, $\delta > 0$ and $dX_t = -\theta\,dt + dB_t$, $\theta > 0$ ($X$ is absorbed when it reaches 0). The local martingale that appears in Itô's formula applied to the process $(e^{-\lambda t}\delta X_t^2)$ is a uniformly integrable martingale whenever $\lambda > 0$. Simple computations show that in the case $\lambda > 0$, problem (2.7) can be reduced to the problem of the form (2.4) with the state space $J = (0,\infty)$ (see the previous remark) and

$$f(x) = (\lambda\delta - 1)x^2 + 2\delta\theta x - \delta.$$

This function $f$ has the form shown in Figure 1 if $\lambda\delta < 1$ and $\lambda\delta + \theta^2\delta > 1$ [cf. Dayanik and Karatzas (2003), Section 6.9, Case III].

Actually, the class of functions $f$ that have the form shown in Figure 1 is a "natural class for which one expects optimal stopping times to be two-sided," though there exist functions $f$ not of this form with corresponding two-sided optimal stopping times and it can happen that optimal stopping times for functions of this form are not two-sided (see Section 3 of the present paper).

2.2. *Free boundary problem and main results.* In order to solve problem (2.4), the free boundary problem (with smooth fit conditions) is usually formulated as follows:

$$\frac{\sigma^2(x)}{2}V''(x) + b(x)V'(x) - \lambda V(x) = -f(x), \qquad x \in (x_1^*, x_2^*); \tag{2.8}$$

$$V(x) = 0, \qquad x \in \mathbb{R} \setminus (x_1^*, x_2^*); \tag{2.9}$$

$$V'_+(x_1^*) = V'_-(x_2^*) = 0, \tag{2.10}$$

where $V'_+$ and $V'_-$ denote, respectively, right and left derivatives of $V$. The form of the free boundary problem (2.8)–(2.10) is motivated by the form of the function $f$. It is natural to expect that the optimal continuation domain here is some interval $(x_1^*, x_2^*)$ [one can also expect that the continuation domain should contain $(x_{1\ell}, x_{2r})$]. The usual way to make use of this free boundary problem is to take an appropriate solution $(V, x_1^*, x_2^*)$ of the problem (2.8)–(2.10) and then to prove that $V = V^*$ and $T_{x_1^*, x_2^*}$ is an optimal stopping time in (2.4).

We say that the free boundary problem (2.8)–(2.10) *loses* a solution of the optimal stopping problem (2.4) if (2.4) has a two-sided optimal stopping time of the form $T_{\alpha,\beta}$ for some real $\alpha < \beta$ and the triplet $(V^*, \alpha, \beta)$ is not a solution of (2.8)–(2.10). This does not happen in many concrete examples with continuous $b$, $\sigma$ and $f$. However, it would be a fairly general situation if $b$, $\sigma$ or $f$ are discontinuous. The reason is that (2.8) is too restrictive in that case: one should not require this equality to be held for all $x \in (x_1^*, x_2^*)$ if one wants to allow discontinuities in $b$, $\sigma$ and $f$.

Below, we shall see that the following modified free boundary formulation is "no-loss" in the sense that it does not lose solutions of (2.4), even if $b$, $\sigma$ and $f$ are allowed to be discontinuous:

$$V' \text{ is absolutely continuous on } [x_1^*, x_2^*]; \tag{2.11}$$

$$\frac{\sigma^2(x)}{2} V''(x) + b(x)V'(x) - \lambda V(x) = -f(x) \quad \text{for } \mu_L\text{-a.a. } x \in (x_1^*, x_2^*); \tag{2.12}$$

$$V(x) = 0, \qquad x \in \mathbb{R} \setminus (x_1^*, x_2^*); \tag{2.13}$$

$$V'_+(x_1^*) = V'_-(x_2^*) = 0. \tag{2.14}$$

We say that a triplet $(V, x_1^*, x_2^*)$ is a *solution* of (2.11)–(2.14) if $x_1^*$ and $x_2^*$ are real numbers, $x_1^* < x_2^*$, $V \in C^1([x_1^*, x_2^*])$ and the triplet $(V, x_1^*, x_2^*)$ satisfies (2.11)–(2.14). Formally, under $V'(x_1^*)$ and $V'(x_2^*)$ in (2.11), one should understand, respectively, right and left derivatives. However, (2.13) and (2.14) imply that the two-sided derivatives exist at both points. We shall see below that condition (2.11) is important.

We say that a solution $(V, x_1^*, x_2^*)$ of (2.11)–(2.14) is *trivial* if $V \equiv 0$. For example, if $x_{1\ell} < x_{1r}$ or $x_{2\ell} < x_{2r}$, then taking any $x_1^* < x_2^*$ belonging either to $[x_{1\ell}, x_{1r}]$ or to $[x_{2\ell}, x_{2r}]$, we get a trivial solution $(0, x_1^*, x_2^*)$. Of course, we are only interested in nontrivial solutions.

The modified free boundary problem (2.11)–(2.14) can be equivalently formulated in the following way, which will be sometimes more convenient for us:

$$V'(x) = \int_{x_1^*}^{x} \frac{2}{\sigma^2(t)} [\lambda V(t) - b(t)V'(t) - f(t)]\, dt, \qquad x \in (x_1^*, x_2^*); \tag{2.15}$$

$$V(x) = 0, \qquad x \in \mathbb{R} \setminus (x_1^*, x_2^*); \tag{2.16}$$

$$V'_-(x_2^*) = 0. \tag{2.17}$$

Similarly, a triplet $(V, x_1^*, x_2^*)$ is a *solution* of (2.15)–(2.17) if $x_1^*$ and $x_2^*$ are real numbers, $x_1^* < x_2^*$, $V \in C^1([x_1^*, x_2^*])$ and the triplet $(V, x_1^*, x_2^*)$ satisfies (2.15)–(2.17). Clearly, a triplet $(V, x_1^*, x_2^*)$ is a solution of (2.11)–(2.14) if and only if it is a solution of (2.15)–(2.17). In connection with (2.15), note that for any function $V \in C^1(\mathbb{R})$, we have $(\lambda V - bV' - f)/\sigma^2 \in L^1_{\mathrm{loc}}(\mathbb{R})$, which follows from (2.2) and (2.3).

The first main result of this paper is a verification theorem for the optimal stopping problem (2.4). Its proof will be given in Section 2.3.

THEOREM 2.1 (Verification theorem). *If $(V, x_1^*, x_2^*)$ is a nontrivial solution of the free boundary problem (2.11)–(2.14), then it is the unique nontrivial solution. $V$ is the value function in the optimal stopping problem (2.4), that is, $V = V^*$ and $T_{x_1^*, x_2^*}$ is the unique optimal stopping time in (2.4).*

REMARKS. (i) Condition (2.11) ensures uniqueness of the nontrivial solution of the free boundary problem (2.11)–(2.14). Problem (2.12)–(2.14) may have nontrivial solutions that have nothing to do with the stopping problem (2.4).

Lamberton and Zervos (2006) prove that value functions of a wide class of stopping problems of the form "maximize $E_x[e^{-A_\tau} g(X_\tau) I(\tau < \infty)]$ over all stopping times $\tau$" are differences of two convex functions. It therefore seems to be of interest whether our free boundary formulation will still have a unique nontrivial solution if we replace (2.11) by the weaker condition "$V$ is the difference of two convex functions." The answer is negative, as the following example shows.

Consider the case $\lambda = 0$ and $b \equiv 0$ and suppose that there exists a nontrivial solution $(V, x_1^*, x_2^*)$ of (2.11)–(2.14) (see Section 3 for necessary and sufficient conditions). We take any continuous function $h : \mathbb{R} \to \mathbb{R}$ such that $h = 0$ on $(-\infty, x_1^*] \cup [x_2^*, \infty)$, $h' = 0$ $\mu_L$-a.e. on $[x_1^*, x_2^*]$, $\int_{x_1^*}^{x_2^*} h(t)\, dt = 0$ and $h$ is not absolutely continuous on $[x_1^*, x_2^*]$ (such a function $h$ can be easily constructed through the Cantor staircase function). We set $H(x) = \int_{-\infty}^x h(t)\, dt$, $x \in \mathbb{R}$ and define the function $\tilde{V}$ by the formula $\tilde{V} = V + H$. Clearly, the triplet $(\tilde{V}, x_1^*, x_2^*)$ satisfies (2.12)–(2.14) and $\tilde{V}$ is the difference of two convex functions. However, $\tilde{V}$ has nothing to do with the stopping problem (2.4) because $\tilde{V} \neq V = V^*$.

(ii) It is interesting to note that we always have strict inequalities $x_1^* < x_{1\ell}$ and $x_2^* > x_{2r}$, regardless of the size of the negative values the function $f$ takes to the left of $x_{1\ell}$ or to the right of $x_{2r}$ (see Proposition 2.9 below). This is different from problems of the form

$$\sup_\tau E_x[e^{-\lambda\tau} g(X_\tau) I(\tau < \infty)].$$

In such problems, a point, where $g$ or $g'$ have a discontinuity, can be a boundary point of the stopping region [for the corresponding examples, see Salminen (1985), page 98, Example (iii), Øksendal and Reikvam (1998), Section 4 or Dayanik and Karatzas (2003), Sections 6.7 and 6.11].

The following result states that the modified free boundary formulation (2.11)–(2.14) is "no-loss" in the sense described above. This justifies the modification of the free boundary as (2.11)–(2.14). In the following theorem, we do not need the structure of the gain function $f$ as specified in Figure 1; we need only the condition $f/\sigma^2 \in L^1_{\mathrm{loc}}(\mathbb{R})$.

THEOREM 2.2. *Let $f$ be any Borel function $\mathbb{R} \to \mathbb{R}$ such that $f/\sigma^2 \in L^1_{\mathrm{loc}}(\mathbb{R})$. If there exist real numbers $x_1^* < x_2^*$ such that $T_{x_1^*,x_2^*}$ is an optimal stopping time in (2.4), then the triplet $(V^*, x_1^*, x_2^*)$ is a solution of (2.11)–(2.14).*

REMARKS. (i) Theorem 2.2 was stated as a conjecture in the first version of this paper. The proof will be given in a subsequent paper (joint with D. Belomestny) which is concerned with this "no-loss" result. The results of Section 3 imply Theorem 2.2 in the particular case $b \equiv 0$ and $\lambda = 0$.

(ii) It follows from Theorems 2.1 and 2.2 that for functions $f$ that have the form shown in Figure 1 and satisfy (2.3), the stopping problem (2.4) has a two-sided optimal stopping time if and only if the free boundary problem (2.11)–(2.14) has a nontrivial solution. It would be interesting to obtain "simple" necessary and sufficient conditions for this in terms of $b$, $\sigma$, $f$ and $\lambda$. We cannot do it in this generality. However, this is done in Section 3 in the particular case $b \equiv 0$ and $\lambda = 0$ (see Theorem 3.1).

(iii) As seen from the discussion above, assumption (2.11) is a key assumption on the value function in the framework of this paper. By Theorem 2.2, (2.11) holds whenever (2.4) has a two-sided optimal stopping time. However, it is interesting to obtain sufficient conditions on the diffusion coefficients, $f$ and $\lambda$ to ensure that the value function $V^*$ in problem (2.4) satisfies (2.11), no matter what form optimal stopping times have.

We cannot solve this problem in general. It follows, however, from the results of Section 3 that in the particular case $b \equiv 0$ and $\lambda = 0$ ($f$ of the form as in Figure 1), the value function $V^*$ satisfies (2.11) if and only if $V^*$ is finite (one can also see necessary and sufficient conditions for this in terms of $\sigma$ and $f$ in Section 3). Compare this with Lamberton and Zervos (2006) and also see Remark (i) after Theorem 2.1.

In the rest of Section 2.2, we study a generalization of our stopping problem and point out some interesting effects. We now consider functions $f$ that have the following form: there exist four segments $I_i = [x_{i,\ell}, x_{i,r}]$, $x_{i,\ell} \le x_{i,r}$, $1 \le i \le 4$, $x_{i,r} < x_{i+1,\ell}$, $1 \le i \le 3$ such that $f = 0$ on $\bigcup_{i=1}^4 I_i$, $f > 0$ on $(x_{1,r}, x_{2,\ell}) \cup (x_{3,r}, x_{4,\ell})$, and $f < 0$ on the rest of real line (we now have two favorable domains compared with one for functions $f$ specified in Figure 1). We suppose that (2.2) and the following modification of (2.3) hold:

$$\frac{f}{\sigma^2} \in L^1_{\mathrm{loc}}(\mathbb{R}), \qquad \frac{1}{f} \in B_{\mathrm{loc}}(x) \qquad \forall x \in \mathbb{R} \setminus \left( \bigcup_{i=1}^{4} I_i \right).$$

To account for the possibility that the stopping region has the form $(-\infty, x_1^*] \cup [x_2^*, x_3^*] \cup [x_4^*, \infty)$, we formulate the modified free boundary problem

(2.18) $V'$ is absolutely continuous on $[x_1^*, x_2^*]$ and on $[x_3^*, x_4^*]$;

(2.19)
$$\frac{\sigma^2(x)}{2} V''(x) + b(x)V'(x) - \lambda V(x) = -f(x)$$
for $\mu_L$-a.a. $x \in (x_1^*, x_2^*) \cup (x_3^*, x_4^*)$;

(2.20) $V(x) = 0, \qquad x \in \mathbb{R} \setminus [(x_1^*, x_2^*) \cup (x_3^*, x_4^*)]$;

(2.21) $V'_+(x_1^*) = V'_-(x_2^*) = V'_+(x_3^*) = V'_-(x_4^*) = 0$

and define its solution $(V, x_1^*, x_2^*, x_3^*, x_4^*)$ in a way similar to that for (2.11)–(2.14).

The following optimal stopping result for the class of functions $f$ with two favorable regions as introduced above can be proven similarly to the proof of Theorem 2.1 (see Section 2.3).

THEOREM 2.3. *If $(V, x_1^*, x_2^*, x_3^*, x_4^*)$ is a solution of the free boundary problem (2.18)–(2.21) such that $V \not\equiv 0$ on $(x_1^*, x_2^*)$ and $V \not\equiv 0$ on $(x_3^*, x_4^*)$, then it is the unique solution with this property. $V$ is the value function in the optimal stopping problem (2.4), that is, $V = V^*$, and the unique optimal stopping time in (2.4) is given by the formula*

$$T_{x_1^*, x_2^*, x_3^*, x_4^*} = \inf\{t \in [0, \infty) : X_t \notin (x_1^*, x_2^*) \cup (x_3^*, x_4^*)\},$$

*where $\inf \varnothing := \infty$. Moreover, $x_1^* < x_{1,\ell}$, $x_2^*, x_3^* \in (x_{2,r}, x_{3,\ell})$ and $x_4^* > x_{4,r}$.*

REMARK. For some functions $f$ of the modified form, the optimal stopping region is two-sided, that is, it is of the form $(-\infty, x_1^*] \cup [x_2^*, \infty)$ [one can think on a function $f$ such that $|f|$ is "small" on $(x_{2,r}, x_{3,\ell})$ and "large" on $\bigcup_{i=1,2,4,5}(x_{i-1,r}, x_{i,\ell})$ with $x_{0,r} := -\infty$ and $x_{5,\ell} := \infty$]. Therefore, it is interesting to understand whether Theorem 2.1 also remains true in this case (to account for two-sided solutions of the stopping problem). The answer is "No"! The form of the function $f$ shown in Figure 1 is crucial for Theorem 2.1.

To illustrate this issue, let us consider some function $\overline{f}$ as in Figure 1 such that the free boundary problem (2.11)–(2.14) with $\overline{f}$ instead of $f$ has a nontrivial solution $(V, x_1^*, x_2^*)$. Then, let us construct a function $f$ of the form we consider now by modifying $\overline{f}$ on $[x_2^*, \infty)$ in such a way that $f > 0$ on some $(x_2^*, x_2^* + \varepsilon)$ [note that $\overline{f} < 0$ on some $(x_2^* - \delta, x_2^*)$ by Remark (ii) after Theorem 2.1]. Clearly, $(V, x_1^*, x_2^*)$ is also a nontrivial solution of (2.11)–(2.14) for this new function $f$, but the stopping time $T_{x_1^*, x_2^*}$ is no more optimal in (2.4) because it is equal to 0 if the starting point belongs to the favorable region $(x_2^*, x_2^* + \varepsilon)$.

2.3. *Auxiliary results and proofs.* Below, we work in the setting of Section 2.1 (in particular, $f$ has the form specified in Figure 1). At first, we need a uniqueness result for the Cauchy problem in (2.15).

Lemma 2.4 (Uniqueness for the Cauchy problem). *Let $I$ be an interval in $\mathbb{R}$ that is either open, semi-open or closed and either bounded or unbounded. Let $g : I \to \mathbb{R}$ be a function such that $g/\sigma^2 \in L^1_{\mathrm{loc}}(I)$. Let $a \in I$, $c \in \mathbb{R}$ and $V$, $\tilde{V}$ be functions $I \to \mathbb{R}$ that satisfy the equation*

$$V'(x) = c + \int_a^x \frac{2}{\sigma^2(t)} [\lambda V(t) - b(t) V'(t) - g(t)]\, dt, \qquad x \in I.$$

*If $V(x_0) = \tilde{V}(x_0)$ and $V'(x_0) = \tilde{V}'(x_0)$ for some $x_0 \in I$, then $V = \tilde{V}$ on $I$.*

Proof. Let us set $U = V - \tilde{V}$ and

$$y = \inf\{x \in I \cap [x_0, \infty) : U(x) \neq 0\}$$

($\inf \varnothing := \infty$) and suppose that $U \not\equiv 0$ on $I \cap [x_0, \infty)$. Then, $y < d$, where $d$ denotes the right endpoint of the interval $I \cap [x_0, \infty)$, and we have

$$U'(x) = \int_y^x \frac{2}{\sigma^2(t)} [\lambda U(t) - b(t) U'(t)]\, dt, \qquad x \in I. \tag{2.22}$$

Let $\delta \in (0, 1]$ be sufficiently small so that $y + \delta < d$ and

$$\max\left( \int_y^{y+\delta} \frac{2\lambda}{\sigma^2(t)}\, dt, \int_y^{y+\delta} \frac{2|b(t)|}{\sigma^2(t)}\, dt \right) \le \frac{1}{3}$$

[see (2.2)], set $m = \sup_{x \in [y, y+\delta]} |U'(x)|$ and note that $m > 0$ and $|U| \le m$ on $[y, y+\delta]$. Now, taking $x \in [y, y+\delta]$ such that $|U'(x)| = m$, we obtain from (2.22) that $m \le \frac{2}{3} m$. This contradiction implies that $U \equiv 0$ on $I \cap [x_0, \infty)$. Similarly, $U \equiv 0$ on $(-\infty, x_0] \cap I$. □

In Lemmas 2.5–2.8 below, we additionally assume that

$$b \text{ is locally bounded on } \mathbb{R}. \tag{2.23}$$

In the following, let $(V, x_1^*, x_2^*)$ be any nontrivial solution of the free boundary problem (2.15)–(2.17).

LEMMA 2.5. *Let $y \in [x_1^*, x_2^*]$ and assume that (2.23) holds.*

(i) *If $V$ attains a local maximum at $y$ and $V(y) \geq 0$, then $f(y) \geq 0$.*
(ii) *If $V$ attains a local minimum at $y$ and $V(y) \leq 0$, then $f(y) \leq 0$.*

PROOF. (i) Since $V'(y) = 0$, we have

$$V'(x) = \int_y^x g(t)\,dt, \qquad x \in [x_1^*, x_2^*],$$

where $g(t) = (2/\sigma^2(t))[\lambda V(t) - b(t)V'(t) - f(t)]$. If $f(y) < 0$, then $g > 0$ in a sufficiently small neighborhood of $y$ [see (2.3) and (2.23)]. Hence, for a sufficiently small $\epsilon > 0$, $V' > 0$ on $(y, y+\epsilon)$ and $V' < 0$ on $(y-\epsilon, y)$ [if $y = x_1^*$ or $y = x_2^*$, one should consider, respectively, only $(y, y+\epsilon)$ or only $(y-\epsilon, y)$]. This contradicts the fact that $V$ attains a local maximum at $y$.

(ii) One can apply the reasoning above to the functions $-V$ and $-f$. □

LEMMA 2.6. *Under assumption (2.23), we have $x_1^* \leq x_{1r}$, $x_2^* \geq x_{2\ell}$ and $V \geq 0$ on $\mathbb{R}$.*

PROOF.

(I) Let us assume that $x_2^* < x_{2\ell}$. There are several cases to consider.
  (1) Suppose that $x_2^* \in (x_{1r}, x_{2\ell})$ and $x_1^* \geq x_{1\ell}$ (see Figure 2). It follows from (2.15)–(2.17), (2.3) and (2.23) that $V < 0$ on $(x_2^* - \epsilon, x_2^*)$ for a sufficiently small $\epsilon > 0$. If $x_1^* < x_{1r}$, then, by Lemma 2.4, $V = 0$ on $[x_1^*, x_{1r}]$. Let us take $y \in [x_1^*, x_2^*]$ such that $V(y) = \inf_{x \in [x_1^*, x_2^*]} V(x)$. The reasoning above ensures that $y \in (x_{1r}, x_{2\ell})$. Since $f > 0$ on $(x_{1r}, x_{2\ell})$, we obtain a contradiction with Lemma 2.5.
  (2) Suppose that $x_2^* \in (x_{1r}, x_{2\ell})$ and $x_1^* < x_{1\ell}$ (see Figure 3). It follows from (2.15)–(2.17), (2.3) and (2.23) that $V < 0$ on $(x_2^* - \epsilon, x_2^*)$ and $V > 0$ on $(x_1^*, x_1^* + \epsilon)$ for a sufficiently small $\epsilon > 0$. By $a$, we denote any point in $(x_1^*, x_2^*)$ such that $V(a) = 0$. Let us take $y \in [x_1^*, a]$ and $z \in [a, x_2^*]$ such that $V(y) = \sup_{x \in [x_1^*, a]} V(x)$ and $V(z) = \inf_{x \in [a, x_2^*]} V(x)$ and note that $y < a < z$ and $V(y) > 0 > V(z)$. By Lemma 2.5, $f(y) \geq 0$ and $f(z) \leq 0$. Due to the form of the function $f$ and the fact that $y < z < x_{2\ell}$, we obtain that $y$, $z \in [x_{1\ell}, x_{1r}]$.

  Since $V'(y) = 0$, we have

$$V'(x) = \int_y^x \frac{2}{\sigma^2(t)}[\lambda V(t) - b(t)V'(t)]\,dt, \qquad x \in [y, z]. \tag{2.24}$$

  In the case where $\lambda > 0$, it follows from (2.24), (2.23) and $V(y) > 0$ that $V' > 0$ on a sufficiently small interval $(y, y+\epsilon)$. This contradicts the fact that $V$ attains a local maximum at $y$.

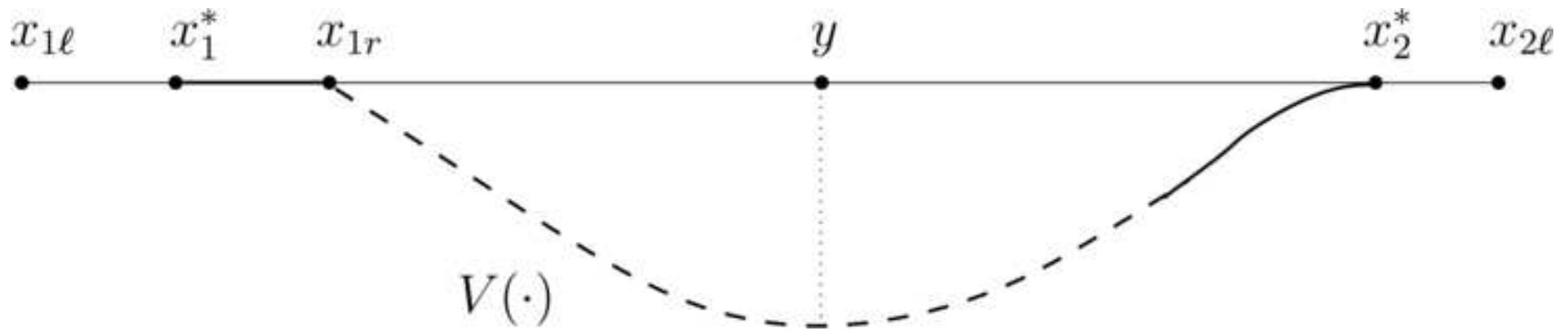

FIG. 2.

Finally, in the case where $\lambda = 0$, the function $\tilde{V}(x) = V(y)$, $x \in [y, z]$ is a solution of (2.24) on the interval $[y, z]$. By Lemma 2.4, $V = \tilde{V}$ on $[y, z]$, hence $V(y) = V(z)$. We obtain a contradiction with $V(y) > 0 > V(z)$.

(3) Suppose that $x_2^* \le x_{1r}$. If $x_1^* \ge x_{1\ell}$, then, by Lemma 2.4, $V \equiv 0$. Since we consider a nontrivial solution $(V, x_1^*, x_2^*)$ of (2.15)–(2.17), we get $x_1^* < x_{1\ell}$. We now obtain a contradiction by reasoning similar to that in part (1). Thus, we have established that $x_2^* \ge x_{2\ell}$. Similarly, $x_1^* \le x_{1r}$.

(II) In the next part, we prove that $V \ge 0$ on $\mathbb{R}$.

Let us first prove that $V \ge 0$ on $[x_{2\ell}, x_2^*]$. If $x_2^* \le x_{2r}$, then, by Lemma 2.4, $V = 0$ on $[x_{2\ell}, x_2^*]$. In the case where $x_2^* > x_{2r}$ (see Figure 4), it follows from (2.15)–(2.17), (2.3) and (2.23) that $V$ is strictly decreasing on a sufficiently small interval $[x_2^* - \epsilon, x_2^*]$, hence $V > 0$ on $[x_2^* - \epsilon, x_2^*)$. To prove that $V$ is decreasing on $[x_{2\ell}, x_2^*]$, assume that this is not the case. There then exist points $y$ and $z$, $x_{2\ell} \le y < z < x_2^*$, such that $V(y) < V(z) = \sup_{x\in[y,x_2^*]} V(x)$. By Lemma 2.5, $f(z) \ge 0$, hence $z \le x_{2r}$. Since $V'(z) = 0$, we obtain from (2.15) that

$$V'(x) = \int_z^x \frac{2}{\sigma^2(t)} [\lambda V(t) - b(t)V'(t)]\, dt, \qquad x \in [y, z]. \tag{2.25}$$

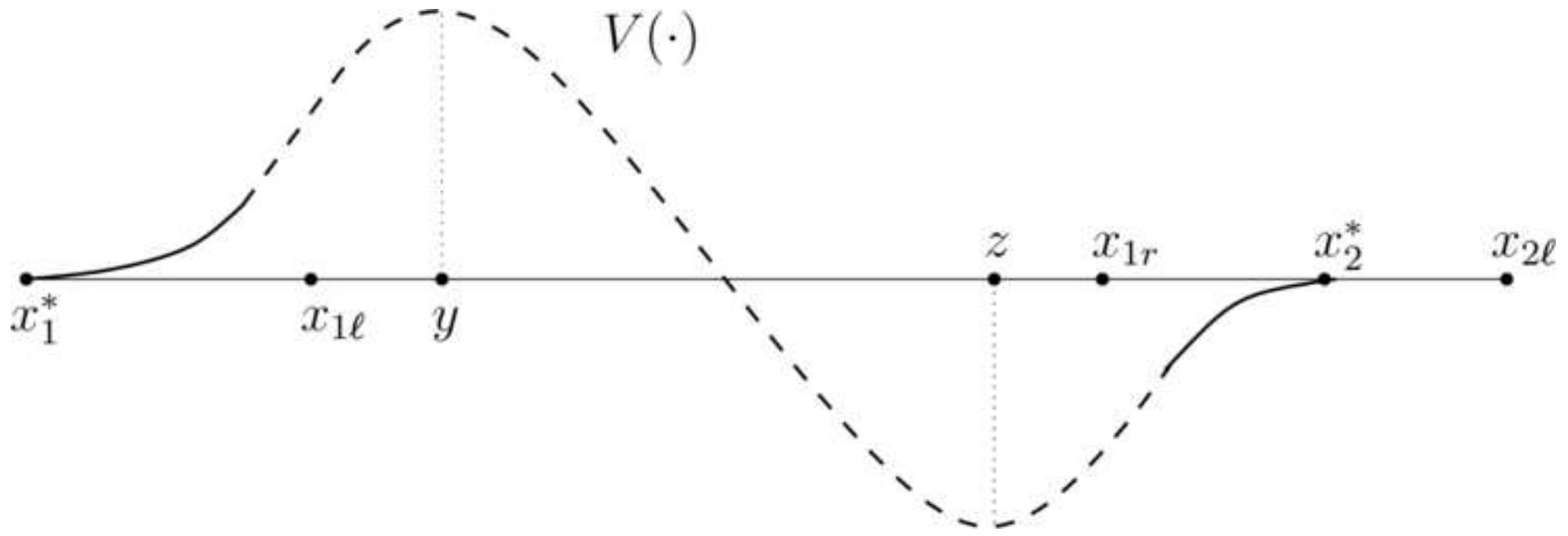

FIG. 3.

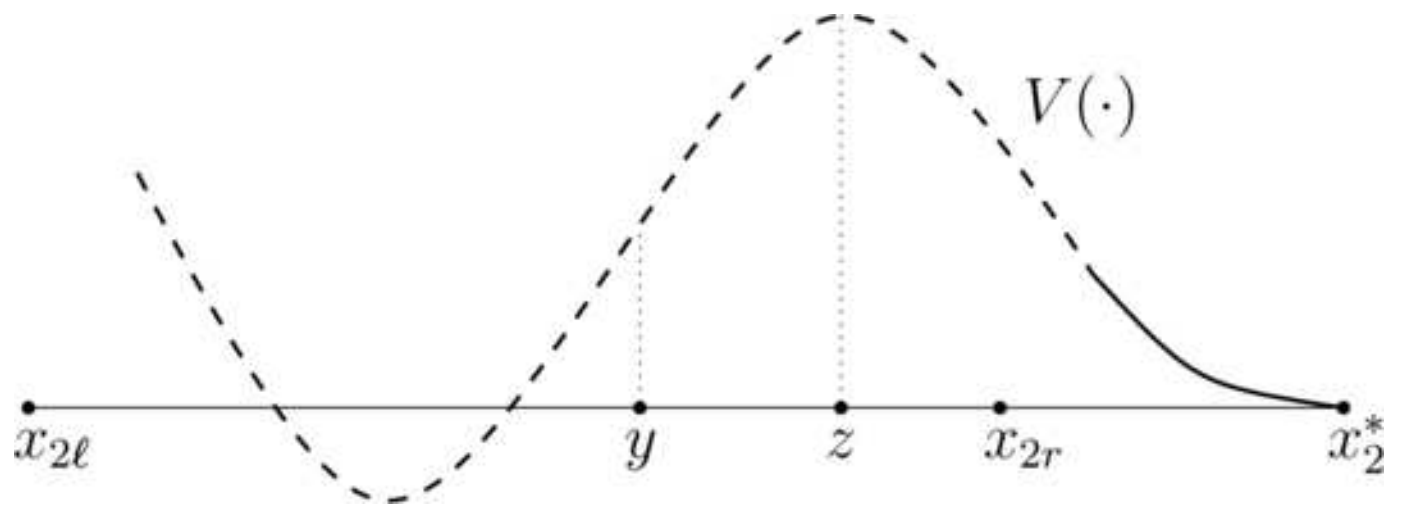

FIG. 4.

If $\lambda > 0$, then $V' < 0$ on some interval $(z-\epsilon, z)$, due to (2.23). This contradicts the fact that $V$ attains a local maximum at $z$. If $\lambda = 0$, then the function $\tilde{V}(x) = V(z)$, $x \in [y, z]$ is a solution of (2.25) on the interval $[y, z]$. By Lemma 2.4, $V = \tilde{V}$ on $[y, z]$, hence $V(y) = V(z)$, and this is a contradiction. Thus, $V \geq 0$ on $[x_{2\ell}, x_2^*]$. Similarly, $V \geq 0$ on $[x_1^*, x_{1r}]$.

Finally, if there exists a point $y \in (x_{1r}, x_{2\ell})$ such that $V(y) < 0$, then we obtain a contradiction with Lemma 2.5 by considering $z \in (x_{1r}, x_{2\ell})$ such that $V(z) = \inf_{x \in [x_{1r}, x_{2\ell}]} V(x)$. This completes the proof. □

In the following lemma, we establish the main part of the verification theorem.

LEMMA 2.7. *Under assumption (2.23), $V^* = V$ and $T_{x_1^*, x_2^*}$ is an optimal stopping time in the stopping problem (2.4).*

PROOF. Let $x \in \mathbb{R}$ be fixed. At first, we prove that we can apply Itô's formula in the standard form to $e^{-\lambda(t \wedge T_{-n,n})} V(X_t^{T_{-n,n}})$, $n \in \mathbb{N}$ (we stop $X$ at $T_{-n,n}$ because it can explode) and obtain

$$\begin{aligned} & e^{-\lambda(t \wedge T_{-n,n})} V(X_{t \wedge T_{-n,n}}) \\ & \qquad = V(x) + \int_0^{t \wedge T_{-n,n}} e^{-\lambda s} V'(X_s) \sigma(X_s)\, dB_s \\ & \qquad\quad + \int_0^{t \wedge T_{-n,n}} e^{-\lambda s} G(X_s)\, ds, \qquad P_x\text{-a.s.}, \ t \in [0, \infty), \end{aligned} \tag{2.26}$$

where $G(y) = \sigma^2(y) V''(y)/2 + b(y) V'(y) - \lambda V(y)$, $y \in \mathbb{R}$. Note that the term $\int_0^{t \wedge T_{-n,n}} e^{-\lambda s} G(X_s)\, ds$ is well defined (though it contains $V''$ that is defined only $\mu_L$-a.e.) because

$$\mu_L(\{t \in [0, \infty) : X_t \in A\}) = 0, \qquad P_x\text{-a.s. whenever } \mu_L(A) = 0, \tag{2.27}$$

which, in turn, can be derived from the occupation times formula [see Revuz and Yor (1999), Chapter VI, Corollary (1.6)].

It follows from (2.11) that $V$ is the difference of two convex functions. By the Itô–Tanaka formula [see Revuz and Yor (1999), Chapter VI, Theorem (1.5)], we get

$$V(X_t^{T_{-n,n}}) = V(x) + \int_0^{t\wedge T_{-n,n}} V'(X_s)\sigma(X_s)\,dB_s$$
$$+ \int_0^{t\wedge T_{-n,n}} V'(X_s)b(X_s)\,ds + \tfrac{1}{2}\int_{\mathbb{R}} L_t^y(X)V'(dy), \tag{2.28}$$
$$P_x\text{-a.s.},\ t\in[0,\infty),$$

where $V'(dy)$ denotes the (signed) Radon measure on $(\mathbb{R},\mathcal{B}(\mathbb{R}))$ with the distribution function $V'$. Since $V'(dy) = V''(y)\,dy$ [see (2.11)], the term with the local time in (2.28) can be rewritten by the occupation times formula as follows:

$$\begin{aligned}\int_{\mathbb{R}} L_t^y(X)V'(dy) &= \int_{\mathbb{R}} L_t^y(X)V''(y)\,dy\\ &= \int_0^t V''(X_s)\,d[X]_s\\ &= \int_0^t V''(X_s)\sigma^2(X_s)\,ds, \qquad P_x\text{-a.s.},\ t\in[0,\infty).\end{aligned}$$

Substituting this in (2.28), we obtain

$$V(X_t^{T_{-n,n}}) = V(x) + \int_0^{t\wedge T_{-n,n}} V'(X_s)\sigma(X_s)\,dB_s$$
$$+ \int_0^{t\wedge T_{-n,n}} [\sigma^2(X_s)V''(X_s)/2 + b(X_s)V'(X_s)]\,ds, \tag{2.29}$$
$$P_x\text{-a.s.},\ t\in[0,\infty).$$

Now, applying Itô's formula to $e^{-\lambda(t\wedge T_{-n,n})}Y_t$, where $Y_t = V(X_t^{T_{-n,n}})$, and using (2.29) we get (2.26). For the sequel, note that in (2.26), we have $G = -fI_{[x_1^*,x_2^*]}$ $\mu_L$-a.e. [see (2.12)].

Since $T_{-n,n}\uparrow\zeta P_x$-a.s. as $n\uparrow\infty$, we have

$$e^{-\lambda(t\wedge T_{-n,n})}V(X_{t\wedge T_{-n,n}}) \to e^{-\lambda t}V(X_t), \qquad P_x\text{-a.s. on } \{\zeta > t\}.$$

Additionally, by our definition $g(\Delta) = 0$ for any function $g:\mathbb{R}\to\mathbb{R}$ and using the fact that $V$ has compact support, we get

$$e^{-\lambda(t\wedge T_{-n,n})}V(X_{t\wedge T_{-n,n}}) \to 0 = e^{-\lambda t}V(X_t), \qquad P_x\text{-a.s. on } \{\zeta \le t\}.$$

Similarly treating the right-hand side of (2.26), we obtain

$$e^{-\lambda t}V(X_t) = V(x) + M_t + \int_0^t e^{-\lambda s}G(X_s)\,ds, \qquad P_x\text{-a.s.},\ t\in[0,\infty), \tag{2.30}$$

where

$$M_t = \int_0^t e^{-\lambda s} V'(X_s)\sigma(X_s)\,dB_s \tag{2.31}$$

is a local martingale. Indeed,

$$\int_0^t [e^{-\lambda s} V'(X_s)\sigma(X_s)]^2\,ds < \infty, \qquad P_x\text{-a.s.},\ t \in [0,\infty)$$

because $\int_0^t \sigma^2(X_s)\,ds < \infty$ $P_x$-a.s. on the set $\{\zeta > t\}$ and, further, $V'$ has compact support.

Let $\tau \in \mathfrak{M}$ be an arbitrary stopping time and let $\tau_n \uparrow \infty$ be a localizing sequence for $M$ (so that each $M^{\tau_n}$ is a uniformly integrable martingale). Setting $F = -G$, we obtain from (2.30) that

$$V(x) = E_x[e^{-\lambda(\tau\wedge\tau_n)} V(X_{\tau\wedge\tau_n})] + E_x \int_0^{\tau\wedge\tau_n} e^{-\lambda s} F(X_s)\,ds. \tag{2.32}$$

(Note that the first term on the right-hand side is finite because $V$ is bounded. Hence, the second term is also finite). By Lebesgue's bounded convergence theorem,

$$E_x[e^{-\lambda(\tau\wedge\tau_n)} V(X_{\tau\wedge\tau_n})] \to E_x[e^{-\lambda\tau} V(X_\tau)], \qquad n \to \infty.$$

Since $F = f I_{[x_1^*, x_2^*]}$ $\mu_L$-a.e., we have $F^+ \le f^+$ $\mu_L$-a.e. and $F^- \le f^-$ $\mu_L$-a.e. Since $\tau \in \mathfrak{M}$, we have

$$E_x \int_0^\tau e^{-\lambda s} F^+(X_s)\,ds < \infty \quad \text{or} \quad E_x \int_0^\tau e^{-\lambda s} F^-(X_s)\,ds < \infty.$$

Separately considering $F^+$ and $F^-$ and applying the monotone convergence theorem, we get

$$E_x \int_0^{\tau\wedge\tau_n} e^{-\lambda s} F(X_s)\,ds \to E_x \int_0^\tau e^{-\lambda s} F(X_s)\,ds, \qquad n \to \infty.$$

Thus, (2.32) implies the definitive equation

$$V(x) = E_x[e^{-\lambda\tau} V(X_\tau)] + E_x \int_0^\tau e^{-\lambda s} F(X_s)\,ds. \tag{2.33}$$

By Lemma 2.6, $V(y) \ge 0$ for all $y \in \mathbb{R}$ and $F(y) \ge f(y)$ for $\mu_L$-a.a. $y \in \mathbb{R}$. Applying (2.27), we obtain from (2.33) that

$$V(x) \ge E_x \int_0^\tau e^{-\lambda s} f(X_s)\,ds$$

for each $\tau \in \mathfrak{M}$. Hence, $V(x) \ge V^*(x)$. Putting $\tau = T_{x_1^*, x_2^*}$ in (2.33), we see that

$$V(x) = E_x \int_0^{T_{x_1^*, x_2^*}} e^{-\lambda s} f(X_s)\,ds.$$

This completes the proof. □

We can now strengthen Lemma 2.6.

LEMMA 2.8. *Under assumption (2.23), we have* $x_1^* \le x_{1\ell}$, $x_2^* \ge x_{2r}$ *and* $V > 0$ *on* $(x_1^*, x_2^*)$.

PROOF. We recall that for the solution $X$ of SDE (2.1) under condition (2.2), it holds that $P_x(T_y < T_z) > 0$ and $P_x(T_z < T_y) > 0$ for any $y < x < z$, where

$$T_y = \inf\{t \in [0, \infty) : X_t = y\}$$

with the usual agreement $\inf \varnothing = \infty$ and where $T_z$ is similarly defined [see Engelbert and Schmidt (1985, 1991) or Karatzas and Shreve (1991), Chapter 5.5.A–B]. Applying Lemma 2.7, we obtain

$$V(x) = V^*(x) \ge E_x \int_0^{T_{x_{1\ell}, x_{2r}}} e^{-\lambda s} f(X_s)\, ds > 0, \qquad x \in (x_{1\ell}, x_{2r}).$$

Hence, $x_1^* \le x_{1\ell}$ and $x_2^* \ge x_{2r}$. Finally, it remains to recall the following fact, which is established in part (II) of the proof of Lemma 2.6: if $x_1^* < x_{1\ell}$ (resp. $x_2^* > x_{2r}$), then $V > 0$ on $(x_1^*, x_{1\ell}]$ (resp. on $[x_{2r}, x_2^*)$). □

PROOF OF THEOREM 2.1. (1) At first, we additionally assume (2.23). We still need to prove the uniqueness of the nontrivial solution of (2.15)–(2.17) and the uniqueness of the optimal stopping time in (2.4).

It follows from Lemmas 2.7 and 2.8 that for any nontrivial solution $(V, x_1^*, x_2^*)$ of (2.15)–(2.17), we have

$$\begin{aligned} V &= V^*, \\ x_1^* &= \sup\{x \le x_{1\ell} : V^*(x) = 0\}, \\ x_2^* &= \inf\{x \ge x_{2r} : V^*(x) = 0\}. \end{aligned}$$

Hence, the nontrivial solution of (2.15)–(2.17) is unique.

Consider any $x \in \mathbb{R}$ and any stopping time $\tau \in \mathfrak{M}$. If $P_x(\tau < T_{x_1^*, x_2^*}) > 0$, then, by Lemma 2.8, $E_x[e^{-\lambda\tau} V(X_\tau)] > 0$. Hence, (2.33) implies that $\tau$ is not optimal in problem (2.4). Now, assume that $P_x(\tau > T_{x_1^*, x_2^*}) > 0$ and consider the process $Y_t = X_{T_{x_1^*, x_2^*} + t} - X_{T_{x_1^*, x_2^*}}$ (note that $T_{x_1^*, x_2^*} < \infty$ $P_x$-a.s.). For any $\epsilon > 0$, we have $\sup_{t \in [0, \epsilon]} Y_t > 0$ $P_x$-a.s. and $\inf_{t \in [0, \epsilon]} Y_t < 0$ $P_x$-a.s. [see, e.g., Karatzas and Shreve (1991), Chapter 5.5.A–B]. It then follows from (2.33), $V \ge 0$ on $\mathbb{R}$ and $F > f$ $\mu_L$-a.e. on $\mathbb{R} \setminus [x_1^*, x_2^*]$ that $\tau$ is not optimal in (2.4). Thus, there exists no other optimal stopping time in (2.4) except $T_{x_1^*, x_2^*}$.

(2) We now prove the result without assuming (2.23). For some fixed $c \in \mathbb{R}$, we consider the scale function of the process $X$

$$p(x) = \int_c^x \exp\left(-\int_c^y \frac{2b(z)}{\sigma^2(z)}\,dz\right)dy, \qquad x \in \mathbb{R}. \tag{2.34}$$

We define the process $\widetilde{X}_t = p(X_t)$, $p(\Delta) := \Delta$, with the state space $J \cup \{\Delta\}$, $J = (p(-\infty), p(\infty))$. We then have

$$d\widetilde{X}_t = \widetilde{\sigma}(\widetilde{X}_t)\,dB_t,$$

with $\widetilde{\sigma}(x) = (p'\sigma) \circ p^{-1}(x)$, $x \in J$. We shall use the alternative notation $\widetilde{P}_x$ for the measure $P_{p^{-1}(x)}$ so that $\widetilde{P}_x(\widetilde{X}_0 = x) = 1$. Consider now the stopping problem

$$\widetilde{V}^*(x) = \sup_{\tau \in \mathfrak{M}} \widetilde{E}_x \int_0^\tau e^{-\lambda s} \widetilde{f}(\widetilde{X}_s)\,ds,$$

where $\widetilde{f} = f \circ p^{-1}$. Clearly, it is a reformulation of the problem (2.4) in the sense that $\widetilde{V}^* = V^* \circ p^{-1}$ and a stopping time $\tau^*$ is optimal in the problem $V^*(x)$ if and only if it is optimal in the problem $\widetilde{V}^*(p(x))$. Note that conditions (2.2) and (2.3) for the functions $\widetilde{b} \equiv 0$, $\widetilde{\sigma}$ and $\widetilde{f}$ are satisfied [one should replace $\mathbb{R}$ with $J$ in (2.2) and (2.3)]. Now, the result follows from part (1) and the fact that the triplet $(V, x_1^*, x_2^*)$ is a nontrivial solution of (2.11)–(2.14) if and only if the triplet $(\widetilde{V}, \widetilde{x}_1^*, \widetilde{x}_2^*) := (V \circ p^{-1}, p(x_1^*), p(x_2^*))$ is a nontrivial solution of the free boundary problem

$$\widetilde{V}' \text{ is absolutely continuous on } [\widetilde{x}_1^*, \widetilde{x}_2^*]; \tag{2.35}$$

$$\frac{\widetilde{\sigma}^2(x)}{2}\widetilde{V}''(x) - \lambda \widetilde{V}(x) = -\widetilde{f}(x) \qquad \text{for } \mu_L\text{-a.a. } x \in (\widetilde{x}_1^*, \widetilde{x}_2^*); \tag{2.36}$$

$$\widetilde{V}(x) = 0, \qquad x \in J \setminus (\widetilde{x}_1^*, \widetilde{x}_2^*); \tag{2.37}$$

$$\widetilde{V}'_+(\widetilde{x}_1^*) = \widetilde{V}'_-(\widetilde{x}_2^*) = 0 \tag{2.38}$$

[we also use Remark (i) at the end of Section 2.1]. □

Finally, we prove the result stated in Remark (ii) after the formulation of Theorem 2.1.

Proposition 2.9. *If $(V, x_1^*, x_2^*)$ is a nontrivial solution of the free boundary problem (2.11)–(2.14), then $x_1^* < x_{1\ell}$ and $x_2^* > x_{2r}$.*

Proof. Since the triplet $(V, x_1^*, x_2^*)$ is a nontrivial solution of (2.11)–(2.14) if and only if the triplet $(\widetilde{V}, \widetilde{x}_1^*, \widetilde{x}_2^*) := (V \circ p^{-1}, p(x_1^*), p(x_2^*))$ is a nontrivial solution of (2.35)–(2.38), we assume, without loss of generality, that $b \equiv 0$.

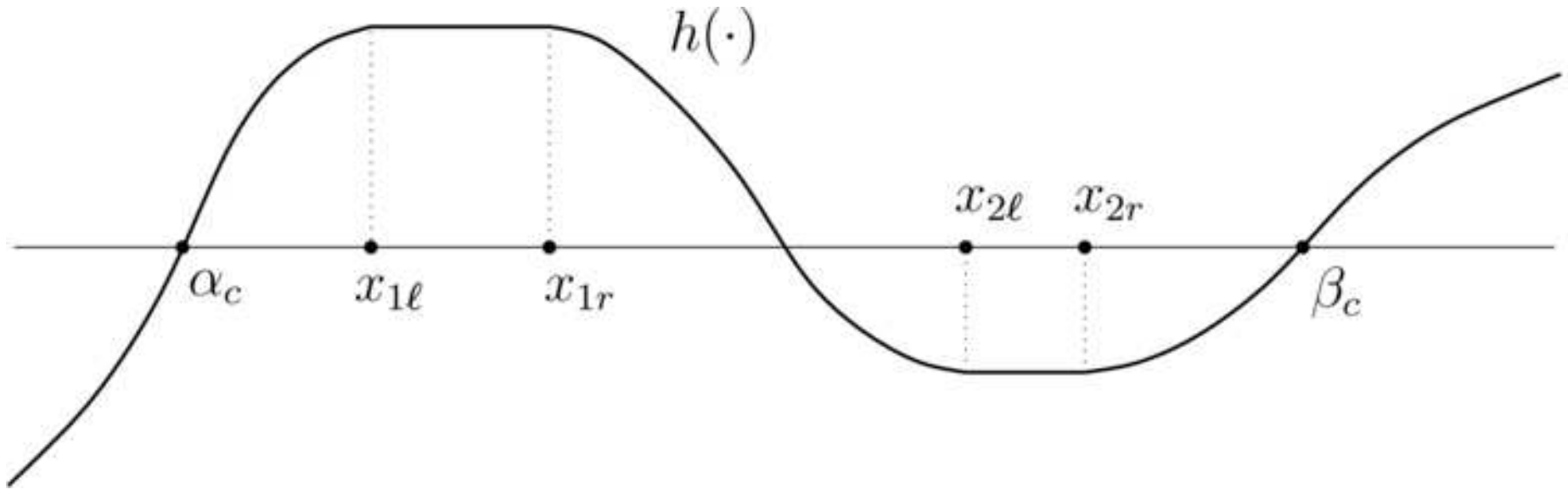

FIG. 5.

By Lemma 2.8, $x_1^* \le x_{1\ell}$. If we suppose that $x_1^* = x_{1\ell}$, then $(V, x_1^*, x_2^*)$ will also be a nontrivial solution of (2.11)–(2.14) with the function $\overline{f}$ instead of $f$, where

$$\overline{f}(x) = \begin{cases} 0, & \text{if } x \in [\overline{x}_{1\ell}, x_{1\ell}], \\ f(x), & \text{otherwise,} \end{cases}$$

with an arbitrary $\overline{x}_{1\ell} < x_{1l}$. Since $\overline{f}$ has the form considered in the paper (see Figure 1), we obtain the contradiction with Lemma 2.8. Thus, $x_1^* < x_{1\ell}$ and, similarly, $x_2^* > x_{2r}$. $\square$

**3. Investigation of the stopping problem in the case $b \equiv 0$ and $\lambda = 0$.** In this section, we consider the case $b \equiv 0$ and $\lambda = 0$ in detail. The assumptions on the functions $f$ and $\sigma$ remain the same as in Section 2. We define the functions $g$ and $h : \mathbb{R} \to \mathbb{R}$ by

$$g(x) = -\frac{2f(x)}{\sigma^2(x)}, \qquad h(x) = \int_0^x g(y)\,dy.$$

The function $h$ is well defined because $f/\sigma^2 \in L^1_{\text{loc}}(\mathbb{R})$ [see (2.3)]. Due to the form of the function $f$, we have $g > 0$ on $(-\infty, x_{1\ell}) \cup (x_{2r}, \infty)$, $g = 0$ on $[x_{1\ell}, x_{1r}] \cup [x_{2\ell}, x_{2r}]$ and $g < 0$ on $(x_{1r}, x_{2\ell})$. Hence, $h$ is strictly increasing on $(-\infty, x_{1\ell}]$ and $[x_{2r}, \infty)$, it is constant on $[x_{1\ell}, x_{1r}]$ and $[x_{2\ell}, x_{2r}]$ and it is strictly decreasing on $[x_{1r}, x_{2\ell}]$. We set $h(\infty) = \lim_{x\to\infty} h(x)$ and $h(-\infty) = \lim_{x\to-\infty} h(x)$.

For any $c \in \mathbb{R}$, we define the function $H(x, c) = h(x) - c$, $x \in \mathbb{R}$. If $c \in \mathbb{R}$ is chosen in such a way that $H(x_{1\ell}, c) > 0$ and $H(-\infty, c) < 0$ [resp. $H(x_{2r}, c) < 0$ and $H(\infty, c) > 0$], then we denote by $\alpha_c$ [resp. $\beta_c$] the unique point in $(-\infty, x_{1\ell})$ [resp. $(x_{2r}, \infty)$] such that $H(\alpha_c, c) = 0$ [resp. $H(\beta_c, c) = 0$]. For an illustration, see Figure 5.

3.1. *Necessary and sufficient conditions for the existence of a nontrivial solution of the free boundary problem* (2.15)–(2.17). We consider the condition $(A_1)$ $h(\infty) > h(-\infty)$ [or, equivalently, $\int_{-\infty}^{\infty} g(y)\,dy > 0$; note that $\int_{-\infty}^{\infty} g(y)\,dy$ is well defined because $\int_{-\infty}^{\infty} g^-(y)\,dy = \int_{x_{1r}}^{x_{2\ell}} g^-(y)\,dy < \infty$].

If $(A_1)$ holds, we additionally introduce the following conditions:

$(A_2)$ If $h(\infty) < h(x_{1\ell})$, then $\int_{\alpha_{h(\infty)}}^{\infty} H(y, h(\infty))\,dy < 0$;

$(A_3)$ If $h(-\infty) > h(x_{2r})$, then $\int_{-\infty}^{\beta_{h(-\infty)}} H(y, h(-\infty))\,dy > 0$.

THEOREM 3.1. *The free boundary problem (2.15)–(2.17) has a nontrivial solution if and only if conditions $(A_1)$–$(A_3)$ hold. In this case, the nontrivial solution is unique.*

PROOF. Assume that conditions $(A_1)$–$(A_3)$ are satisfied. We set $m_1 = h(x_{2r}) \vee h(-\infty)$ and $m_2 = h(x_{1\ell}) \wedge h(\infty)$. There exist $c_1$ and $c_2$, $m_1 < c_1 < c_2 < m_2$, such that

$$\int_{\alpha_{c_1}}^{\beta_{c_1}} H(y, c_1)\,dy > 0 \quad \text{and} \quad \int_{\alpha_{c_2}}^{\beta_{c_2}} H(y, c_2)\,dy < 0.$$

There then exists $c^* \in (c_1, c_2)$ such that

$$\int_{\alpha_{c^*}}^{\beta_{c^*}} H(y, c^*)\,dy = 0.$$

It is now clear that the triplet $(V, \alpha_{c^*}, \beta_{c^*})$ is a nontrivial solution of (2.15)–(2.17), where

$$V(x) = \begin{cases} \displaystyle\int_{\alpha_{c^*}}^{x} H(y, c^*)\,dy, & \text{if } x \in (\alpha_{c^*}, \beta_{c^*}), \\ 0 & \text{otherwise.} \end{cases} \tag{3.1}$$

The converse and the uniqueness can also be easily verified (alternatively, the uniqueness follows from Theorem 2.1). □

It is a remarkable fact that the value function of the optimal stopping problem can be determined in explicit form [see (3.1)] based on the free boundary formulation (2.15)–(2.17). This shows the usefulness of the modified formulation of the free boundary problem.

3.2. *Study of the optimal stopping problem when $(A_1)$–$(A_3)$ are not satisfied.* Suppose that at least one of the conditions $(A_1)$–$(A_3)$ is violated. It is also interesting to consider the stopping problem (2.4) in this case. For the sequel, note that our assumption $b \equiv 0$ implies that the solution $X$ of (2.1) does not explode and is recurrent [see Engelbert and Schmidt (1985, 1991) or Karatzas and Shreve (1991), Chapter 5.5.A]. If at least one of conditions $(A_1)$–$(A_3)$ is violated, then we are in the situation of exactly one of the following three cases:

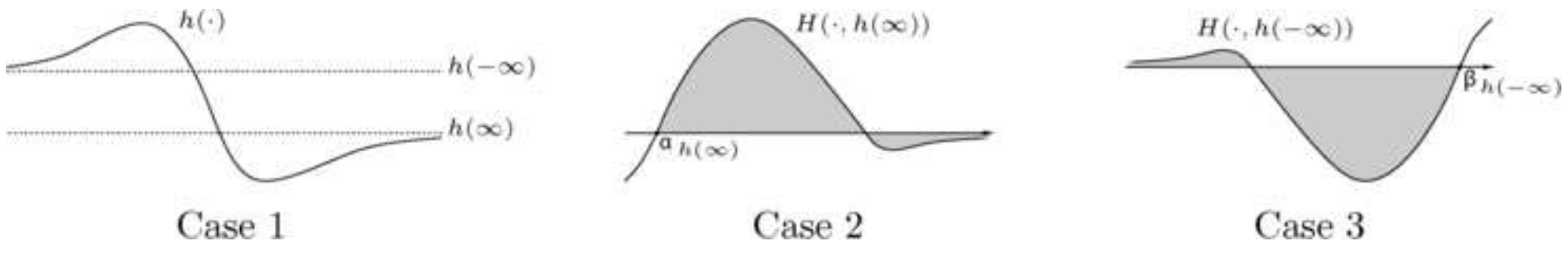

Fig. 6.

*Case* 1: $h(\infty) \le h(-\infty)$ [or, equivalently, $\int_{-\infty}^{\infty} g(y)\,dy \le 0$];
*Case* 2: $h(-\infty) < h(\infty) < h(x_{1\ell})$ and $\int_{\alpha_{h(\infty)}}^{\infty} H(y, h(\infty))\,dy \ge 0$;
*Case* 3: $h(x_{2r}) < h(-\infty) < h(\infty)$ and $\int_{-\infty}^{\beta_{h(-\infty)}} H(y, h(-\infty))\,dy \le 0$.

Note that each of these cases excludes the other ones.

For an illustration, see Figure 6.

3.3. *Study of case* 1. In case 1, $h(\infty)$ and $h(-\infty)$ are finite. We set

$$m = \frac{h(\infty) + h(-\infty)}{2}, \qquad K^+ = \int_{\mathbb{R}} H^+(y, m)\,dy, \qquad K^- = \int_{\mathbb{R}} H^-(y, m)\,dy.$$

Recall that for any real numbers $\alpha < \beta$, $T_{\alpha,\beta} \in \mathfrak{M}$ and, moreover,

$$E_x \int_0^{T_{\alpha,\beta}} |f(X_s)|\,ds < \infty, \qquad x \in \mathbb{R}$$

(see Lemma A.3).

Theorem 3.2 (Solution of the stopping problem in case 1).

(i) *For any $\tau \in \mathfrak{M}$ and $x \in \mathbb{R}$, we have*

$$E_x \int_0^{\tau} f(X_s)\,ds < V^*(x),$$

*that is, there exists no optimal stopping time.*

(ii) *There exist sequences $a_n \downarrow -\infty$ and $b_n \uparrow \infty$ such that for any $x \in \mathbb{R}$,*

$$E_x \int_0^{T_{a_n,b_n}} f(X_s)\,ds \to V^*(x), \tag{3.2}$$

*that is, the sequence of stopping times $\{T_{a_n,b_n}\}$ is asymptotically optimal.*

(iii) *If $K^+ = \infty$ or $K^- = \infty$, then*

$$V^*(x) = \infty, \qquad x \in \mathbb{R}.$$

*If $K^- \le K^+ < \infty$, then*

$$V^*(x) = \int_{-\infty}^{x} H(y, m)\,dy, \qquad x \in \mathbb{R}.$$

*If $K^+ \le K^- < \infty$, then*

$$V^*(x) = -\int_{x}^{\infty} H(y, m)\,dy, \qquad x \in \mathbb{R}.$$

In particular, we have

$$K^+ \vee K^- < \infty \iff V^*(x) < \infty \qquad \forall x \in \mathbb{R},$$
$$K^+ \vee K^- = \infty \iff V^*(x) = \infty \qquad \forall x \in \mathbb{R}.$$

Let us also note that $K^+ \vee K^-$ can be finite only if $h(\infty) = h(-\infty)$.

Remarks. (i) The situation of case 1 can be heuristically interpreted as the situation when the "negative tails" of the function $f$ are light compared with the "positive midst." This interpretation suggests that it is never optimal to stop (because $X$ is recurrent and $\lambda = 0$), which agrees with Theorem 3.2. In this connection, we note that $\tau \equiv \infty$ is not an optimal stopping time here because $\tau \notin \mathfrak{M}$ (see Lemma 3.3).

(ii) We would like to remark that the asymptotic optimality of $T_{a_n,b_n}$ in (3.2) is not true for all sequences $a_n \downarrow -\infty$ and $b_n \uparrow \infty$. At the end of this subsection, we present a corresponding example.

The proof of Theorem 3.2 will follow from Lemmas 3.6–3.8 below. At first, however, we need several auxiliary results.

Lemma 3.3. *For the stopping time $\tau \equiv \infty$, we have $\tau \notin \mathfrak{M}$.*

Proof. Since $X$ is a recurrent continuous local martingale, we have $[X]_\infty = \infty$, $P_x$-a.s. Then, for the local time of $X$, we have $L^y_\infty = \infty$ $P_x$-a.s. for any $y \in \mathbb{R}$ [Revuz and Yor (1999), Chapter VI, Example (1.27)]. By the occupation times formula,

$$\int_0^\infty f^+(X_s)\,ds = \int_0^\infty \frac{f^+(X_s)}{\sigma^2(X_s)}\,d[X]_s$$
$$= \int_{\mathbb{R}} \frac{f^+(y)}{\sigma^2(y)} L^y_\infty(X)\,dy = \infty, \qquad P_x\text{-a.s.}$$

Similarly, $\int_0^\infty f^-(X_s)\,ds = \infty$, $P_x$-a.s. $\square$

In Lemmas 3.4–3.8 below, we assume that

$$K^+ \ge K^-. \tag{3.3}$$

The case $K^+ \le K^-$ can be treated similarly.

Lemma 3.4. *Let $z \in \mathbb{R}$ be an arbitrary real number. There exist sequences $\{a_n\}$, $\{b_n\}$ and $\{c_n\}$ such that the following statements hold:*

(i) $a_n \downarrow -\infty$, $a_1 \le x_{1\ell}$;
(ii) $b_n \uparrow \infty$, $b_1 \ge x_{2r}$;

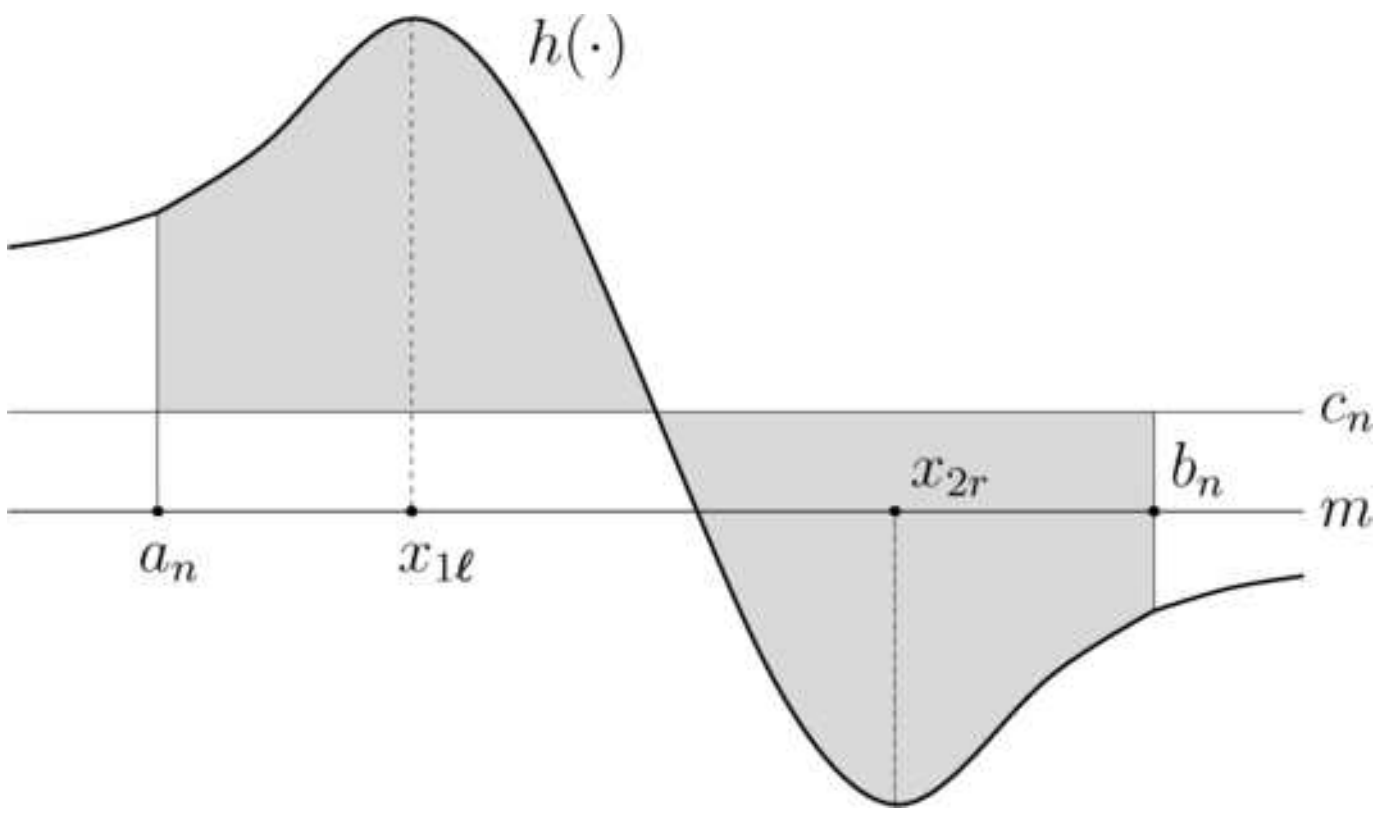

FIG. 7.

(iii) $m \le c_n < h(a_n)$;
(iv) $\int_{a_n}^{b_n} H(y, c_n)\,dy = 0$;
(v) *if* $a_n < z$, *then* $\int_{a_n}^{z} H(y, c_n)\,dy \ge \int_{a_n}^{z} H(y, m)\,dy - \frac{1}{n}$.

In connection with statement (v), let us note that $a_n < z$ for sufficiently large $n$. For an illustration, see Figure 7.

PROOF OF LEMMA 3.4. At first, we take any sequences $a'_n \downarrow -\infty$ and $b'_n \uparrow \infty$ such that $a'_1 \le x_{1\ell}$, $b'_1 \ge x_{2r}$ and, for any $n$,

$$\int_{a'_n}^{b'_n} H(y, m)\,dy \ge 0.$$

This can be done due to (3.3). We now construct the sequence $\{c'_n\}$ in the following way. If $\int_{a'_n}^{b'_n} H(y, m)\,dy = 0$, we take $c'_n = m$. If $\int_{a'_n}^{b'_n} H(y, m)\,dy > 0$, we take $c'_n$ sufficiently close to $m$ so that the following properties are satisfied:

(a) $m < c'_n < h(a'_n)$;
(b) $\int_{a'_n}^{b'_n} H(y, c'_n)\,dy \ge 0$;
(c) if $a'_n < z$, then $\int_{a'_n}^{z} H(y, c'_n)\,dy \ge \int_{a'_n}^{z} H(y, m)\,dy - \frac{1}{n}$.

Since $c'_n > m \ge h(\infty)$, we have $\int_{\mathbb{R}} H^-(y, c'_n)\,dy = \infty$. Consequently, there exists $b''_n \ge b'_n$ such that

$$\int_{a'_n}^{b''_n} H(y, c'_n)\,dy = 0.$$

Finally, let us denote by $\{b_n\}$ any monotone subsequence of $\{b''_n\}$ and by $\{a_n\}$ and $\{b_n\}$ the corresponding subsequences of $\{a'_n\}$ and $\{c'_n\}$. Clearly, statements (i)–(v) hold. $\square$

Now, let $\{a_n\}$, $\{b_n\}$ and $\{c_n\}$ be any sequences satisfying conditions (i)–(v) of Lemma 3.4. We consider the optimal stopping problem

$$V_n^*(x) = \sup_{\tau \le T_{a_n,b_n}} E_x \int_0^\tau f(X_s)\,ds, \tag{3.4}$$

where the supremum is taken over all stopping times $\tau \le T_{a_n,b_n}$ (note that by Lemma A.3, $\tau \in \mathfrak{M}$ whenever $\tau \le T_{a_n,b_n}$). We define

$$V_n(x) = \begin{cases} \displaystyle\int_{a_n}^x H(y, c_n)\,dy, & \text{if } x \in (a_n, b_n),\\ 0, & \text{otherwise.} \end{cases}$$

$V_n$ is then continuous and $V_n > 0$ on $(a_n, b_n)$ (see Figure 7).

LEMMA 3.5. *$V_n$ is identical to the optimal stopping value in (3.4), that is, $V_n = V_n^*$, and $T_{a_n,b_n}$ is the unique optimal stopping time.*

PROOF. If $x \notin (a_n, b_n)$, then the statement is clear; thus, let $x \in (a_n, b_n)$. Let $\widetilde{V}_n$ be any function such that $\widetilde{V}_n \in C^1(\mathbb{R}) \cap C^2((-\infty, a_n] \cup [b_n, \infty))$ and $\widetilde{V}_n = V_n$ on $[a_n, b_n]$. Note that we cannot take $\widetilde{V}_n = V_n$ because $V_n'$ has discontinuities at the points $a_n$ and $b_n$. We can apply Itô's formula in the standard form to $\widetilde{V}_n(X)$ (as in the proof of Lemma 2.7) and obtain

$$\begin{aligned}\widetilde{V}_n(X_t) = \widetilde{V}_n(x) &+ \int_0^t \widetilde{V}_n'(X_s)\sigma(X_s)\,dB_s\\ &+ \tfrac{1}{2}\int_0^t \widetilde{V}_n''(X_s)\sigma^2(X_s)\,ds, \qquad P_x\text{-a.s.},\ t \in [0,\infty).\end{aligned}$$

Hence, by (2.12),

$$V_n(X_t) = V_n(x) + M_t - \int_0^t f(X_s)\,ds, \qquad P_x\text{-a.s. on } \{t \le T_{a_n,b_n}\}, \tag{3.5}$$

where

$$M_t = \int_0^t V_n'(X_s)\sigma(X_s)\,dB_s$$

(we used also (2.27) because $V''(y)\sigma^2(y)/2 = -f(y)$ only for $\mu_L$-a.a. $y \in [a_n, b_n]$). By Lemma A.1 and boundedness of $V_n'$ on $(a_n, b_n)$, $M^{T_{a_n,b_n}}$ is a uniformly integrable martingale. Hence, we obtain from (3.5) that for any stopping time $\tau \le T_{a_n,b_n}$,

$$V_n(x) = E_x V_n(X_\tau) + E_x \int_0^\tau f(X_s)\,ds.$$

This implies the statement of Lemma 3.5. □

LEMMA 3.6. *For any* $x \in \mathbb{R}$,

$$E_x \int_0^{T_{a_n,b_n}} f(X_s)\,ds = V_n^*(x) \uparrow V^*(x), \qquad n \uparrow \infty.$$

PROOF. The equality is a part of Lemma 3.5. The sequence $\{V_n^*(x)\}_{n\in\mathbb{N}}$ is increasing and for each $n$, $V_n^*(x) \le V^*(x)$. By the monotone convergence theorem applied separately to $f^+$ and $f^-$, we have, for any $\tau \in \mathfrak{M}$,

$$E_x \int_0^{\tau} f(X_s)\,ds = \lim_{n\to\infty} E_x \int_0^{\tau\wedge T_{a_n,b_n}} f(X_s)\,ds \le \lim_{n\to\infty} V_n^*(x).$$

Hence, $V_n^*(x) \uparrow V^*(x)$. $\square$

Let us recall that we still assume that $K^+ \ge K^-$ [see (3.3)].

LEMMA 3.7. *If* $K^+ = \infty$, *then*

$$V^*(x) = \infty, \qquad x \in \mathbb{R}. \tag{3.6}$$

*If* $K^+ < \infty$, *then*

$$V^*(x) = \int_{-\infty}^{x} H(y,m)\,dy, \qquad x \in \mathbb{R}. \tag{3.7}$$

PROOF. Let us recall that

$$\int_{a_n}^{z} H(y,m)\,dy - \frac{1}{n} \le \int_{a_n}^{z} H(y,c_n)\,dy \le \int_{a_n}^{z} H(y,m)\,dy$$

(see Lemma 3.4). By Lemmas 3.5 and 3.6,

$$\begin{aligned} V^*(z) &= \lim_{n\to\infty} V_n^*(z) = \lim_{n\to\infty} V_n(z) = \lim_{n\to\infty} \int_{a_n}^{z} H(y,c_n)\,dy \\ &= \lim_{n\to\infty} \int_{a_n}^{z} H(y,m)\,dy = \int_{-\infty}^{z} H(y,m)\,dy. \end{aligned}$$

Since $z$ is an arbitrary point in $\mathbb{R}$, we obtain (3.7). Finally, it remains to note that the right-hand side of (3.7) is identically infinite if $K^+ = \infty$. $\square$

LEMMA 3.8. *For any* $\tau \in \mathfrak{M}$ *and* $x \in \mathbb{R}$,

$$E_x \int_0^{\tau} f(X_s)\,ds < V^*(x).$$

PROOF. We prove this result by contradiction. Assume that there exists $\tau \in \mathfrak{M}$ and $x \in \mathbb{R}$ such that

$$E_x \int_0^{\tau} f(X_s)\,ds = V^*(x). \tag{3.8}$$

(1) At first, we consider the case $V^*(x) < \infty$. By Lemma 3.3, there exists a sufficiently large $n$ such that $P_x(\tau < T_{a_n,b_n}) > 0$. By Lemma A.3, $\tau \vee T_{a_n,b_n} \in \mathfrak{M}$. To obtain a contradiction, it is enough to prove that

$$E_x \int_0^{\tau \vee T_{a_n,b_n}} f(X_s)\,ds > E_x \int_0^{\tau} f(X_s)\,ds. \tag{3.9}$$

If $N$ is a process such that $N^{T_{a_n,b_n}}$ is a uniformly integrable martingale, then

$$E_x[N_\tau I(\tau < T_{a_n,b_n})] = E_x[N_{T_{a_n,b_n}} I(\tau < T_{a_n,b_n})]. \tag{3.10}$$

We rewrite (3.5) as

$$V_n(x) + M_t = V_n(X_t) + \int_0^t f(X_s)\,ds \qquad P_x\text{-a.s. on } \{t \le T_{a_n,b_n}\}$$

and substitute the process $N_t = V_n(x) + M_t$ in (3.10). We obtain

$$E_x\left[\left(V_n(X_\tau) + \int_0^\tau f(X_s)\,ds\right) I(\tau < T_{a_n,b_n})\right] \tag{3.11}$$

$$= E_x\left[\int_0^{T_{a_n,b_n}} f(X_s)\,ds I(\tau < T_{a_n,b_n})\right]. \tag{3.12}$$

Since we assume (3.8) and $V^*(x) < \infty$, $E_x \int_0^\tau f(X_s)\,ds$ is finite. Hence, $E_x[\int_0^\tau f(X_s)\,ds I(\tau < T_{a_n,b_n})]$ is finite. Together with (3.11) and $E_x[V_n(X_\tau) \times I(\tau < T_{a_n,b_n})] > 0$, this implies that

$$E_x\left[\int_0^{T_{a_n,b_n}} f(X_s)\,ds I(\tau < T_{a_n,b_n})\right] > E_x\left[\int_0^{\tau} f(X_s)\,ds I(\tau < T_{a_n,b_n})\right]. \tag{3.13}$$

Finally, we add the finite quantity $E_x[\int_0^\tau f(X_s)\,ds I(\tau \ge T_{a_n,b_n})]$ to both sides of (3.13) and obtain (3.9).

(2) We now consider the case $V^*(x) = \infty$. By the occupation times formula, we have

$$E_x \int_0^\tau f^+(X_s)\,ds = \int_{x_{1r}}^{x_{2\ell}} \frac{f^+(y)}{\sigma^2(y)} E_x L_\tau^y(X)\,dy \tag{3.14}$$

[note that $f^+ = 0$ outside $(x_{1r}, x_{2\ell})$] and, similarly,

$$E_x \int_0^\tau f^-(X_s)\,ds = \int_{\mathbb{R}} \frac{f^-(y)}{\sigma^2(y)} E_x L_\tau^y(X)\,dy. \tag{3.15}$$

It follows from (3.8) that $E_x \int_0^\tau f^+(X_s)\,ds = \infty$. Lemma A.2, formula (3.14) and local integrability of $f^+/\sigma^2$ imply that $E_x L_\tau^y(X) = \infty$, $\forall y \in \mathbb{R}$. Hence, (3.15) yields $E_x \int_0^\tau f^-(X_s)\,ds = \infty$. This contradicts $\tau \in \mathfrak{M}$. $\square$

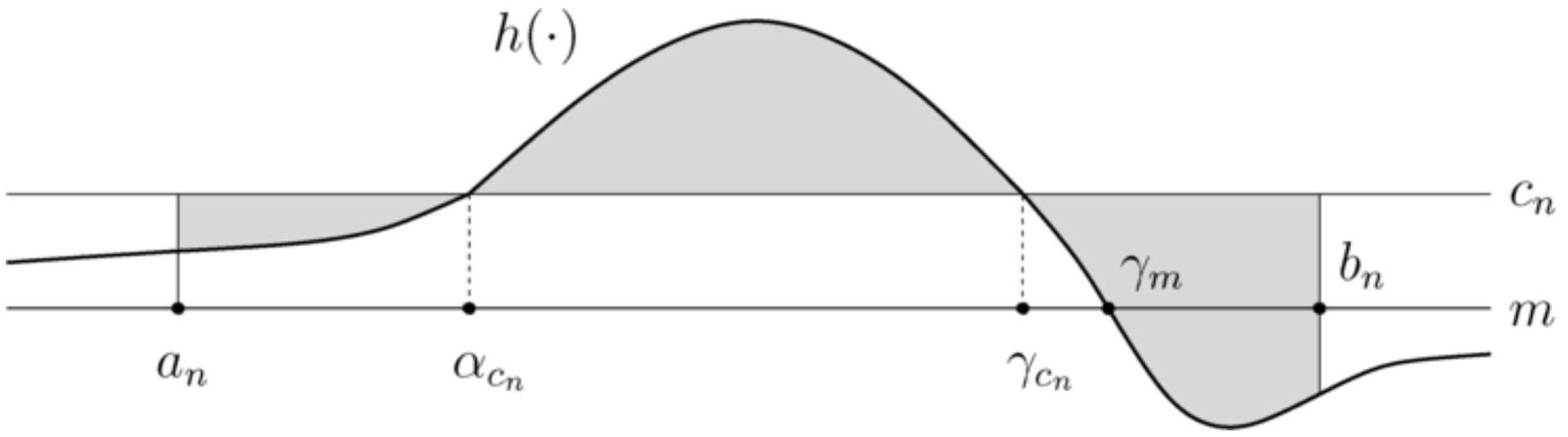

FIG. 8.

This concludes the proof of Theorem 3.2 in the case $K^+ \geq K^-$. The case $K^+ \leq K^-$ can be dealt with in a similar way. We omit the details.

EXAMPLE. Consider the situation $h(\infty) = h(-\infty)$, $K^+ = \infty$ and $K^- < \infty$. We shall construct sequences $a_n \downarrow -\infty$ and $b_n \uparrow \infty$ such that

$$\overline{\lim_{n\to\infty}} E_x \int_0^{T_{a_n,b_n}} f(X_s)\, ds < V^*(x) \tag{3.16}$$

for all $x \in \mathbb{R}$ [of statement (ii) of Theorem 3.2].

The construction is illustrated by Figure 8. For $c \in (h(x_{2\ell}), h(x_{1r}))$, we use the notation $\gamma_c$ for the unique point in $(x_{1r}, x_{2\ell})$ such that $h(\gamma_c) = c$. Let $\{b_n\}$ be a sequence, $b_n \uparrow \infty$ and $b_1 \geq x_{2r}$. Further, let $\{c_n\}$ be a sequence such that $c_n \downarrow m$ and:

(a) $m < c_n < h(x_{1r})$;
(b) $\int_{\alpha_{c_n}}^{\gamma_{c_n}} H(y, c_n)\, dy > K^-$ [note that $\int_{\alpha_c}^{\gamma_c} H(y, c)\, dy \uparrow K^+ = \infty$ as $c \downarrow m$];
(c) $\int_{\gamma_{c_n}}^{b_n} (-H(y, c_n))\, dy \leq K^-$ [note that $\int_{\gamma_m}^{b_n} (-H(y, m))\, dy < K^-$ for all $n$].

For each $n$, we have $\int_{\alpha_{c_n}}^{b_n} H(y, c_n)\, dy > 0$. Since $\int_{-\infty}^{\alpha_{c_n}} H(y, c_n)\, dy = -\infty$, there exists $a_n < \alpha_{c_n}$ such that

$$\int_{a_n}^{b_n} H(y, c_n)\, dy = 0.$$

Clearly, $a_n \to -\infty$. Without loss of generality, we assume that the sequence $\{a_n\}$ is monotone (otherwise, we consider a monotone subsequence) and define

$$U_n(x) = \begin{cases} \displaystyle\int_{a_n}^{x} H(y, c_n)\, dy, & \text{if } x \in (a_n, b_n), \\ 0, & \text{otherwise.} \end{cases}$$

Note that $U_n$ is continuous, strictly decreasing on $[a_n, \alpha_{c_n}]$ and $[\gamma_{c_n}, b_n]$ and strictly increasing on $[\alpha_{c_n}, \gamma_{c_n}]$. In particular,

$$\sup_{x\in\mathbb{R}} U_n(x) = U_n(\gamma_{c_n}) \leq K^-, \tag{3.17}$$

by condition (c) above. One can verify, with the help of Itô's formula, that

$$E_x \int_0^{T_{a_n,b_n}} f(X_s)\,ds = U_n(x), \qquad x \in \mathbb{R}. \tag{3.18}$$

Since $K^+ = \infty$, we have $V^* \equiv \infty$. Together with (3.17) and (3.18), this implies (3.16).

3.4. *Study of cases* 2 *and* 3. We recall that in case 2, $h(-\infty) < h(\infty) < h(x_{1\ell})$ and $\int_{\alpha_{h(\infty)}}^{\infty} H(y, h(\infty))\,dy \ge 0$, and in case 3, $h(x_{2r}) < h(-\infty) < h(\infty)$ and $\int_{-\infty}^{\beta_{h(-\infty)}} H(y, h(-\infty))\,dy \le 0$. For real numbers $\alpha$, $\beta$, we define the one-sided stopping times

$$T_\alpha^- = \inf\{t \in [0,\infty) : X_t \le \alpha\}, \tag{3.19}$$

$$T_\beta^+ = \inf\{t \in [0,\infty) : X_t \ge \beta\}, \tag{3.20}$$

(as usual, $\inf \varnothing = \infty$). It is important that $T_\alpha^-, T_\beta^+ \in \mathfrak{M}$ (see Lemma A.4). We introduce the functions

$$V^-(x) = \begin{cases} 0, & \text{if } x \le \alpha_{h(\infty)}, \\ \displaystyle\int_{\alpha_{h(\infty)}}^{x} H(y, h(\infty))\,dy, & \text{if } x > \alpha_{h(\infty)}, \end{cases}$$

and

$$V^+(x) = \begin{cases} \displaystyle -\int_{x}^{\beta_{h(-\infty)}} H(y, h(-\infty))\,dy & \text{if } x < \beta_{h(-\infty)}, \\ 0, & \text{if } x \ge \beta_{h(-\infty)}. \end{cases}$$

THEOREM 3.9 (Solution of the stopping problem in cases 2 and 3). *In cases* 2 *and* 3 *the optimal stopping value* $V^*$ *is given by* $V^* = V^-$ *(resp.* $V^* = V^+$*). The optimal stopping times are unique and are given by the one-sided stopping times* $T^-_{\alpha_{h(\infty)}}$ *(resp.* $T^+_{\beta_{h(-\infty)}}$*).*

The proof of Theorem 3.9 is omitted. It can be obtained similarly to that of Lemma 2.7 and Theorem 2.1. The details of that proof concerning a possible explosion of $X$ can be omitted here. Note that, unlike case 1, the value function $V^*$ is finite in cases 2 and 3.

REMARK. Case 2 can be heuristically interpreted as the situation when the "right negative tail" of the function $f$ is light, while the "left negative tail" of $f$ is heavy. This interpretation makes it natural that the optimal stopping time should have the form $T_\alpha^-$ for a suitably chosen $\alpha$. The situation in case 3 is symmetric.

## APPENDIX

Here, we prove some technical lemmas which are used in the proofs and which also seem to be of independent interest.

Below, $J=(\ell,r)$, $-\infty\le\ell<r\le\infty$, and $X$ is a (possibly explosive) $J\cup\{\Delta\}$-valued diffusion that satisfies the SDE (2.1) and starts at the point $x\in J$ under the measure $P_x$ ($X$ explodes when it tends either to $\ell$ or to $r$ at a finite time). The coefficients $b$ and $\sigma$ are Borel functions $J\to\mathbb{R}$ that satisfy

$$\sigma(x)\neq 0\ \forall x\in J \qquad \frac{1}{\sigma^2}\in L^1_{\mathrm{loc}}(J), \qquad \frac{b}{\sigma^2}\in L^1_{\mathrm{loc}}(J), \tag{A.1}$$

where $L^1_{\mathrm{loc}}(J)$ denotes the class of functions $J\to\mathbb{R}$ that are integrable on compact subintervals of $J$. Let us define the strictly increasing function $p$ by formula (2.34) and the process $\widetilde{X}_t=p(X_t)$, $p(\Delta):=\Delta$, with the state space $\widetilde{J}\cup\{\Delta\}$, $\widetilde{J}=(\widetilde{\ell},\widetilde{r}):=(p(\ell),p(r))$. We then have

$$d\widetilde{X}_t=\widetilde{\sigma}(\widetilde{X}_t)\,dB_t,$$

with $\widetilde{\sigma}(x)=(p'\sigma)\circ p^{-1}(x)$, $x\in\widetilde{J}$. Note that condition (A.1) with $\widetilde{J}$ instead of $J$ holds for the functions $\widetilde{b}\equiv 0$ and $\widetilde{\sigma}$. We shall use the alternative notation $\widetilde{P}_x$ for the measure $P_{p^{-1}(x)}$ so that $\widetilde{P}_x(\widetilde{X}_0=x)=1$. For $\alpha<\beta$ in $J$, we use the notation $T_{\alpha,\beta}$ of (2.6). For $\alpha<\beta$ in $\widetilde{J}$, we define

$$\widetilde{T}_{\alpha,\beta}:=\inf\{t\in[0,\infty):\widetilde{X}_t\le\alpha \text{ or } \widetilde{X}_t\ge\beta\} \qquad (=T_{p^{-1}(\alpha),p^{-1}(\beta)}).$$

Lemma A.1. *For any $p>0$, $\alpha,\beta\in J$, $\alpha<\beta$, we have*

$$E_x\biggl(\int_0^{T_{\alpha,\beta}}\sigma^2(X_s)\,ds\biggr)^p<\infty, \qquad x\in J \tag{A.2}$$

*(or, equivalently, $E_x[X]^p_{T_{\alpha,\beta}}<\infty$, $x\in J$).*

Proof. If $x\notin(\alpha,\beta)$, then the statement is clear. Let us assume that $x\in(\alpha,\beta)$ and set $\widetilde{\varkappa}(y)=\sigma\circ p^{-1}(y)$, $y\in\widetilde{J}$. Below, we denote positive constants used in estimates by $c_1$, $c_2$, etc. Note that $\widetilde{\varkappa}(y)\le c_1\widetilde{\sigma}(y)$, $y\in[p(\alpha),p(\beta)]$. We have

$$\begin{aligned} E_x\biggl(\int_0^{T_{\alpha,\beta}}\sigma^2(X_s)\,ds\biggr)^p &= \widetilde{E}_{p(x)}\biggl(\int_0^{\widetilde{T}_{p(\alpha),p(\beta)}}\widetilde{\varkappa}^2(\widetilde{X}_s)\,ds\biggr)^p \\ &\le c_2\widetilde{E}_{p(x)}\biggl(\int_0^{\widetilde{T}_{p(\alpha),p(\beta)}}\widetilde{\sigma}^2(\widetilde{X}_s)\,ds\biggr)^p. \end{aligned}$$

Hence, it is enough to prove (A.2) under the additional assumption $b \equiv 0$. Then, $(X_{t\wedge T_{\alpha,\beta}})$ is a bounded martingale. For $q \geq 1$, Burkholder–Davis–Gundy inequalities yield

$$E_x[X]^{q/2}_{T_{\alpha,\beta}} \leq c_3 E_x \left(\sup_{t\geq 0} X_{t\wedge T_{\alpha,\beta}}\right)^q < \infty.$$

This completes the proof. $\square$

Below, $L^y_t(X)$ denotes the local time of $X$ at time $t$ and level $y$.

Lemma A.2. *Let $x \in J$ and $\tau$ be an arbitrary stopping time. Consider the function $h : J \to [0,\infty]$ defined by $h(y) = E_x L^y_\tau(X)$. Then, either $h(y) = \infty$ $\forall y \in J$ or $h$ is bounded.*

Let us stress that neither finiteness nor boundedness of $\tau$ is assumed.

Proof of Lemma A.2. By Revuz and Yor (1999), Chapter VI, Exercise (1.23),

$$L^{p(y)}_\tau(p(X)) = p'(y)L^y_\tau(X), \qquad P_x\text{-a.s.}, \ y \in J.$$

Hence,

$$h(y) = \frac{1}{p'(y)}\widetilde{E}_{p(x)}L^{p(y)}_\tau(\widetilde{X}), \qquad y \in J.$$

Therefore, it is enough to prove the lemma under the additional assumption $b \equiv 0$.

For some sequences $a_n \downarrow \ell$ and $b_n \uparrow r$, set $h_n(y) = E_x L^y_{\tau\wedge T_{a_n,b_n}}(X)$, $y \in J$. Since the local time remains unchanged after the explosion time, $h_n(y) \uparrow h(y)$. Assume that $h(y_0) < \infty$ for some $y_0 \in J$ and consider an arbitrary $y \in J$. By the Tanaka formula under the measure $P_x$ [see Revuz and Yor (1999), Chapter VI, Theorem (1.2)],

$$\begin{aligned}|X_{\tau\wedge T_{a_n,b_n}} - y| = |x-y| &+ \int_0^{\tau\wedge T_{a_n,b_n}} \operatorname{sgn}(X_s - y)\sigma(X_s)\,dB_s \\ &+ L^y_{\tau\wedge T_{a_n,b_n}}(X), \qquad P_x\text{-a.s.},\end{aligned} \tag{A.3}$$

where

$$\operatorname{sgn} y = \begin{cases} 1, & \text{if } y > 0, \\ -1, & \text{if } y \leq 0. \end{cases}$$

For each $y$, the process $M_t = \int_0^{t\wedge T_{a_n,b_n}} \operatorname{sgn}(X_s - y)\sigma(X_s)\,dB_s$ is a uniformly integrable martingale, by Lemma A.1. Taking the expectation in (A.3), we get

$$E_x|X_{\tau\wedge T_{a_n,b_n}} - y| = |x-y| + h_n(y). \tag{A.4}$$

In particular,

$$E_x|X_{\tau\wedge T_{a_n,b_n}} - y_0| = |x - y_0| + h_n(y_0). \tag{A.5}$$

Since $h(y_0) < \infty$, we obtain from (A.5) that $c := \sup_n E_x|X_{\tau\wedge T_{a_n,b_n}}| < \infty$. Now, (A.4) implies that, for any $n$,

$$h_n(y) = E_x[|X_{\tau\wedge T_{a_n,b_n}} - y| - |x - y|] \le E_x|X_{\tau\wedge T_{a_n,b_n}} - x| \le c + |x|.$$

Hence, the function $h$ is bounded. □

LEMMA A.3. *For any Borel function $f : J \to \mathbb{R}$ such that $f/\sigma^2 \in L^1_{\mathrm{loc}}(J)$ and any $\alpha, \beta \in J$, $\alpha < \beta$, we have*

$$E_x \int_0^{T_{\alpha,\beta}} |f(X_s)|\, ds < \infty, \qquad x \in J.$$

PROOF. We need only to consider the case $x \in (\alpha, \beta)$. Using the occupation times formula (under the measure $P_x$), we obtain

$$\int_0^{T_{\alpha,\beta}} |f(X_s)|\, ds = \int_0^{T_{\alpha,\beta}} \frac{|f(X_s)|}{\sigma^2(X_s)}\, d[X]_s = \int_\alpha^\beta \frac{|f(y)|}{\sigma^2(y)} L^y_{T_{\alpha,\beta}}(X)\, dy. \tag{A.6}$$

By Lemma A.2, the function $y \mapsto E_x L^y_{T_{\alpha,\beta}}(X)$ is bounded on $J$ (that is because this function equals 0 at the point $y = \alpha$). Since we have $f/\sigma^2 \in L^1_{\mathrm{loc}}(J)$, the statement of the lemma follows from (A.6). □

Below, we use the notation $T^-_\alpha$ and $T^+_\beta$, $\alpha, \beta \in J$, for one-sided stopping times, as in (3.19) and (3.20).

LEMMA A.4. *Let $\alpha, \beta \in J$. For any Borel function $f : J \to \mathbb{R}$ that has the form shown in Figure 1 and which satisfies $f/\sigma^2 \in L^1_{\mathrm{loc}}(J)$, we have*

$$E_x \int_0^{T^-_\alpha} f^+(X_s)\, ds < \infty \quad \text{and} \quad E_x \int_0^{T^+_\beta} f^+(X_s)\, ds < \infty, \qquad x \in J.$$

PROOF. The proof is similar to that of Lemma A.3. The form of $f$ being as shown in Figure 1 ensures that the integral on the right-hand side of the analogue of (A.6) can be taken over the compact subinterval $[x_{1\ell}, x_{2r}]$ of $J$. □

**Acknowledgments.** We are grateful to H. R. Lerche for interesting discussions and for drawing our attention to the relevant literature. We express deep thanks to M. Zervos for an essential idea for the proof of Theorem 2.2 and to P. Bank and D. Lamberton for interesting discussions which gave rise to removal of the assumption of local boundedness of $b$ that was present in the first version of this paper. We thank the anonymous referee for several relevant comments that helped to improve the paper.

## REFERENCES

Beibel, M. and Lerche, H. R. (2000). Optimal stopping of regular diffusions under random discounting. *Theory Probab. Appl.* **45** 547–557. MR1968720

Bensoussan, A. and Lions, J.-L. (1973). Problèmes de temps d'arrêt optimal et inéquations variationnelles paraboliques. *Appl. Anal.* **3** 267–294. MR0449843

Dayanik, S. (2003). Optimal stopping of linear diffusions with random discounting. Working paper. Available at http://www.princeton.edu/~sdayanik/papers/additive.pdf.

Dayanik, S. and Karatzas, I. (2003). On the optimal stopping problem for one-dimensional diffusions. *Stochastic Process. Appl.* **107** 173–212. MR1999788

Engelbert, H. J. and Schmidt, W. (1985). On one-dimensional stochastic differential equations with generalized drift. *Stochastic Differential Systems (Marseille-Luminy, 1984). Lecture Notes in Control and Inform. Sci.* **69** 143–155. Springer, Berlin. MR0798317

Engelbert, H. J. and Schmidt, W. (1991). Strong Markov continuous local martingales and solutions of one-dimensional stochastic differential equations. III. *Math. Nachr.* **151** 149–197. MR1121203

Friedman, A. (1976). *Stochastic Differential Equations and Applications*. **2**. Academic Press, New York. MR0494491

Glowinski, R., Lions, J.-L. and Trémolières, R. (1976). *Analyse numérique des inéquations variationnelles. Tome 2*. Dunod, Paris.

Graversen, S. E., Peskir, G. and Shiryaev, A. N. (2000). Stopping Brownian motion without anticipation as close as possible to its ultimate maximum. *Teor. Veroyatnost. i Primenen.* **45** 125–136. MR1810977

Karatzas, I. and Ocone, D. (2002). A leavable bounded-velocity stochastic control problem. *Stochastic Process. Appl.* **99** 31–51. MR1894250

Karatzas, I. and Shreve, S. E. (1991). *Brownian Motion and Stochastic Calculus*, 2nd ed. Springer, New York. MR1121940

Lamberton, D. and Zervos, M. (2006). On the problem of optimally stopping a one-dimensional Itô diffusion. Unpublished manuscript.

Nagai, H. (1978). On an optimal stopping problem and a variational inequality. *J. Math. Soc. Japan* **30** 303–312. MR0488321

Øksendal, B. and Reikvam, K. (1998). Viscosity solutions of optimal stopping problems. *Stochastics Stochastics Rep.* **62** 285–301. MR1613260

Peskir, G. and Shiryaev, A. N. (2006). *Optimal Stopping and Free-Boundary Problems*. Birkhäuser, Basel. MR2256030

Revuz, D. and Yor, M. (1999). *Continuous Martingales and Brownian Motion*, 3rd ed. Springer, Berlin. MR1725357

Salminen, P. (1985). Optimal stopping of one-dimensional diffusions. *Math. Nachr.* **124** 85–101. MR0827892

Zabczyk, J. (1984). Stopping games for symmetric Markov processes. *Probab. Math. Statist.* **4** 185–196. MR0792784

Zhang, X. (1994). Analyse numérique des options américaines dans un modèle de diffusion avec sauts. Ph.D. thesis, CERMA–Ecole Nationale des Ponts et Chaussées.

Department of Mathematical Stochastics
University of Freiburg
Eckerstrasse 1
79104 Freiburg
Germany
E-mail: ruschen@stochastik.uni-freiburg.de

Institute of Mathematics
MA 7–4
Berlin University of Technology
Strasse des 17. Juni 136
10623 Berlin
and
Quantitative Products Laboratory
Global Markets Equity
Deutsche Bank AG
Alexanderstrasse 5
10178 Berlin
Germany
E-mail: urusov@math.tu-berlin.de